\documentclass[10pt]{article}
\usepackage{amsmath,amssymb}
\usepackage{url}

\newcommand{\nwc}{\newcommand}
\nwc{\draftdate}{\today\quad PLEASE DO NOT CIRCULATE}
\newtheorem{theorem}{Theorem}
\newtheorem{lemma}{Lemma}

\newtheorem{corollary}{Corollary}

\nwc{\qref}[1]{(\ref{#1})} %\eqref seems to have a spacing problem

\renewcommand{\div}{\nabla\!\cdot\!}
\nwc{\udotgrad}{\uu\!\cdot\!\nabla}

\nwc{\dsl}{\displaystyle}

\newcommand{\PP}{{\cal P}}
\newcommand{\QQ}{{\cal Q}}

\nwc{\pcon}{\beta}  % The constant formerly known as \gamma
\nwc{\pavs}{\bar{p}_s}
\nwc{\epss}{{\varepsilon}}
\nwc{\pgh}{p_{gh}}
\nwc{\SP}{{\cal S}_p}
\nwc{\SPb}{{\cal S}_\Gamma}
\nwc{\cut}{\xi}

\nwc{\Gdag}{G^{\dag}}
\nwc{\I}{{\cal I}}
\nwc{\gam}{\gamma}
\nwc{\yy}{\vec{y}}
\nwc{\zz}{\vec{z}}
\nwc{\hu}{\hat{u}}
\nwc{\hv}{\hat{v}}

\nwc{\xuh}{X_h}
\nwc{\xph}{Y_h}
\nwc{\Tmax}{T_{\rm max}}
\nwc{\Phifac}{\Phi\grad\nn}
\nwc{\tf}{\tilde{f}}
\nwc{\lapg}{\Delta_\G}
\nwc{\gradg}{\nabla_\G}
\nwc{\lapG}{\Delta_\Gamma}
\nwc{\lapI}{\Delta_\I}
\nwc{\gradG}{\nabla_\Gamma}
\nwc{\gradI}{\nabla_\I}
\nwc{\N}{\mathbb{N}}
\nwc{\ip}[1]{\langle #1 \rangle}
\nwc{\Omgs}{{\Omega_s}}
\nwc{\G}{{\cal G}}

\newcommand{\uuD}{\uu_{\Delta t}}
\newcommand{\uud}{\vec{U}_{\Delta t}}
\nwc{\ffd}{\ff_{\Dt}}
\nwc{\UU}{\vec{U}}
\nwc{\emb}{\hookrightarrow}
\nwc{\shalf}{{\textstyle\frac12}}

\nwc{\uspace}{{H^2\cap H^1_0(\Omega,\R^N)}}

\nwc{\uuin}{\uu_{\rm in}}
\nwc{\hin}{h_{\rm in}}

\nwc{\Htwoloc}{H^2_{\rm loc}}

\newcommand{\Hdiv}{H({\rm div};\Omega)}
\newcommand{\Hdivs}{H({\rm div};\Omega_s)}
\newcommand{\Hcurl}{H({\rm curl};\Omega)}
\newcommand{\bau}{{\cal A}\uu}
\newcommand{\dist}{\mathop{\rm dist}\nolimits}
\newcommand{\lap}{\Delta }
\newcommand{\Omg}{\Omega}
\newcommand{\bOmg}{\bar{\Omg}}
\newcommand{\Gam}{\Gamma}

\newcommand{\eps}{\varepsilon}

\nwc{\ba}{\vec{a}}
\newcommand{\nn}{\vec{n}}

\newcommand{\ff}{\vec{f}}
\newcommand{\tff}{\tilde{f}}
\newcommand{\tuu}{\tilde{u}}
\renewcommand{\gg}{\vec{g}}
\newcommand{\uu}{\vec{u}}
\newcommand{\vv}{\vec{v}}
\newcommand{\ww}{\vec{w}}
\newcommand{\bb}{\vec{b}}

\newcommand{\xx}{\vec{x}}

\newcommand{\curl}{\nabla\times}
\newcommand{\nproj}{\nn\nn^t}
\newcommand{\tanproj}{I-\nproj}
\newcommand{\upe}{\vec{u}_{\perp}}
\newcommand{\upa}{\vec{u}_{\parallel}}

\newcommand{\grad}{\nabla}
\newcommand{\nder}{\nn\cdot\!\grad}
\newcommand{\pa}{\partial}
\newcommand{\D}{\partial}

\newcommand{\R}{\mathbb{R}}
\newcommand{\C}{\mathbb{C}}

\renewcommand{\P}{\PP}

\newcommand{\Dt}{\Delta t}
\newcommand{\noind}{\noindent}
\newcommand{\pe}{p_{\mbox{\tiny E}}}
\newcommand{\ps}{p_{\mbox{\tiny S}}}
\newcommand{\bps}{\bar{p}_{\mbox{\tiny S}}}
\newcommand{\qe}{q_{\mbox{\tiny E}}}
\newcommand{\qs}{q_{\mbox{\tiny S}}}

\renewcommand{\div}{\nabla \cdot}
\newcommand{\ndot}{\nn \cdot}
\renewcommand{\(}{\big(}
\renewcommand{\)}{\big)}
\newcommand{\<}{\big<}
\renewcommand{\>}{\big>}

\begin{document}

\title{Divorcing pressure from viscosity in\\
incompressible Navier-Stokes dynamics}
\author{Jian-Guo~Liu\textsuperscript{1}, Jie~Liu\textsuperscript{2}, and
Robert~L.~Pego\textsuperscript{3}}
\maketitle

{\abstract{
We show that in bounded domains with no-slip boundary conditions,
the Navier-Stokes pressure can be determined in a such way that it is
strictly dominated by viscosity. As a consequence, in a general domain
we can treat the Navier-Stokes equations as a perturbed vector diffusion
equation, instead of as a perturbed Stokes system.  We illustrate the
advantages of this view in a number of ways. In particular, we provide
simple proofs of (i) local-in-time existence and uniqueness of strong
solutions for an unconstrained formulation of the Navier-Stokes
equations, and (ii) the unconditional stability and convergence of
difference schemes that are implicit only in viscosity and explicit in
both pressure and convection terms, requiring no solution of
stationary Stokes systems or inf-sup conditions.  }}

%\bigskip \begin{center} PLEASE DO NOT CIRCULATE \end{center}

\bigskip
%Key words and phrases:

\footnotetext[1]{Department of Mathematics \&
Institute for Physical Science and Technology,
University of Maryland, College Park MD 20742.
Email: jliu@math.umd.edu}
\footnotetext[2]{Department of Mathematics,
University of Maryland, College Park MD 20742.
Email: jieliu@math.umd.edu}
\footnotetext[3]
{Department of Mathematical Sciences, Carnegie Mellon University,
Pittsburgh, PA 15213.
Email: rpego@cmu.edu}

%\input{nse-organization}
%\pagebreak

%%%-------------------------------------
\section{Introduction} \label{S.introduction}

The pressure term has always created problems for understanding the
Navier-Stokes equations of incompressible flow.  Pressure plays a role
like a Lagrange multiplier to enforce the incompressibility
constraint, and this has been a main source of difficulties.
Our general aim in this paper is to show that the pressure can be
obtained in a way that leads to considerable simplifications in both
computation and analysis.

From the computational point of view, typical difficulties are related to
the lack of an evolution equation for updating the pressure
dynamically and the lack of useful boundary conditions for
determining the pressure by solving boundary-value problems.
Existing methods able to handle these difficulties are
sophisticated and lack the robustness and flexibility that would
be useful to address more complex problems. For example, finite element
methods have required carefully arranging approximation spaces for
velocity and pressure to satisfy an inf-sup compatibility
condition \cite{GR}.  Projection methods too have typically
encountered problems related to low-accuracy approximation of
the pressure near boundaries \cite{Ch,Te2,OID}.
Yet much of the scientific and technological significance of the
Navier-Stokes equations derives from their role in the modeling of
physical phenomena such as lift, drag, boundary-layer separation and
vortex shedding, for which the behavior of the pressure near
boundaries is of great importance.

Our main results in this article indicate that in bounded domains
with no-slip boundary conditions, the Navier-Stokes pressure can be
determined in a such way that it is {\em strictly dominated by viscosity}.
To explain, let us take $\Omg$ to be a bounded, connected
domain in $\R^N$ ($N\ge2$) with $C^3$ boundary $\Gamma=\pa\Omg$.
The Navier-Stokes equations for incompressible fluid flow in
$\Omg$ with no-slip boundary conditions on $\Gamma$ take the form
\begin{align}
  \pa_t \uu + \udotgrad \uu + \nabla p &= \nu \Delta \uu + \ff
  \quad\mbox{  in $\Omg$},  \label{NSE1} \\
  \div \uu &= 0
  \quad\mbox{  in $\Omg$}, \label{NSE2} \\
  \quad \uu &=0  \quad\mbox{ on $\Gamma$}.  \label{NSE3}
\end{align}
Here $\uu$ is the fluid velocity, $p$ the pressure, and $\nu$ is
the kinematic viscosity coefficient, assumed to be a fixed
positive constant.

A standard way to determine $p$ is via the Helmholtz-Hodge
decomposition. We let $\PP$ denote the Helmholtz projection operator onto
divergence-free fields, and recall that it is defined as follows.
Given any $\ba \in L^2(\Omg, \R^N)$, there is a
unique $q \in H^1(\Omg)$ with $\int_{\Omg} q = 0$ such that
$\P\ba :=\ba+\grad q$ satisfies
\begin{equation}\label{hodge}
 0= \int_\Omg (\P\ba)\cdot\grad\phi =
 \int_\Omg (\ba+\grad q)\cdot\grad\phi 
\quad\mbox{for all $\phi\in H^1(\Omg)$.}
\end{equation}
The pressure $p$ in
\qref{NSE1} is determined by taking $\ba= \udotgrad \uu - \ff -
\nu \Delta \uu$. Then \qref{NSE1} is rewritten as
\begin{equation}
 \pa_t \uu + \PP( \udotgrad \uu - \ff - \nu \Delta \uu) = 0. \label{NSE4}
\end{equation}

In this formulation, solutions formally satisfy $\pa_t(\div\uu)=0$.
Consequently the zero-divergence condition \qref{NSE2} needs to be
imposed only on initial data. 
Nevertheless, the pressure is determined from \qref{NSE4} in principle
even for velocity fields that do not respect the incompressibility constraint.
However, the dissipation in \qref{NSE4} appears degenerate due to the
fact that $\PP$ annihilates gradients, so the analysis of \qref{NSE4}
is usually restricted to spaces of divergence-free fields.

Alternatives are possible in which the
pressure is determined differently when the velocity field has
non-zero divergence. Instead of \qref{NSE4}, we propose to consider
\begin{equation}
 \pa_t \uu + \PP( \udotgrad \uu - \ff - \nu \Delta \uu )
 = \nu \nabla(\div \uu). \label{NSE5}
\end{equation}
Of course there is no difference as long as $\div\uu=0$.
But we argue that \qref{NSE5} enjoys superior stability properties,
for two reasons.  The first is heuristic. The
incompressibility constraint is enforced in a more robust way, because
the divergence of velocity satisfies a weak form of the diffusion
equation with no-flux (Neumann) boundary conditions ---
Due to \qref{hodge}, for all appropriate test functions $\phi$ we have
\begin{equation}
  \int_\Omg \pa_t \uu \cdot\nabla \phi =
  \nu \int_\Omg \nabla(\div\uu) \cdot \nabla \phi.
  \label{heat1}
\end{equation}
Taking $\phi=\div\uu$ we get the dissipation identity
\begin{equation}
  \frac{d}{dt} \frac12 \int_\Omg (\div\uu)^2  +
 \nu \int_\Omg |\nabla(\div\uu)|^2 = 0.
  \label{heat2}
\end{equation}
Due to the Poincar\'e inequality and the fact that $\int_\Omg\div\uu=0$,
the divergence of velocity is smoothed and decays exponentially in $L^2$ norm.
Naturally, if $\div\uu = 0$ initially, this remains true for all
later time, and one has a solution of the standard Navier-Stokes
equations \qref{NSE1}--\qref{NSE3}.

The second reason is much deeper. To explain, we recast
\qref{NSE5} in the form \qref{NSE1} while explicitly identifying
the separate contributions to the pressure term made by the
convection and viscosity terms. Using the Helmholtz projection
operator $\PP$, we introduce the {\em Euler pressure} $\pe$ and
{\em Stokes pressure} $\ps$ via the relations
%\begin{equation}
\begin{align}
& \PP( \udotgrad \uu - \ff) = \udotgrad \uu - \ff +
\nabla \pe , \label{EulerP}
\\[4pt] %\end{equation} %\begin{equation}
& \PP(-\Delta \uu ) = -\Delta \uu
+\nabla(\div\uu) + \nabla \ps . \label{StokesP}
\end{align}
%\end{equation}
This puts \qref{NSE5} into the form \qref{NSE1} with $p=\pe+\nu\ps$:
\begin{equation}
 \pa_t \uu + \udotgrad \uu + \nabla \pe + \nu \nabla \ps
 = \nu \Delta \uu + \ff.  \label{NSE6}
\end{equation}
Identifying the Euler and Stokes pressure terms in this way allows one
to focus separately on the difficulties peculiar to each.
The Euler pressure is nonlinear, but of lower order.
Since the Helmholtz projection is orthogonal, naturally the Stokes
pressure satisfies
\begin{equation}
\int_\Omg |\grad\ps|^2 \le \int_\Omg |\lap\uu|^2 
\qquad\mbox{if $\div\uu=0$.}
\end{equation}
The key observation is that the Stokes pressure term
is actually {\em strictly} dominated by the viscosity term,
regardless of the divergence constraint.
We regard the following theorem as the main achievement of this paper.
\begin{theorem} \label{T.main}
Let $\Omg \subset \R^N$ ($N\ge2$) be a connected bounded domain with $C^3$
boundary. Then for any $\eps>0$, there exists $C\ge 0$
such that for all vector fields $\uu \in \uspace$,
the Stokes pressure $\ps$ determined by
\qref{StokesP} satisfies
\begin{equation}
 \int_\Omg |\nabla \ps |^2 \le %\left(\frac23+\eps\right) 
 \pcon \int_\Omg |\lap
   \uu|^2 + C \int_\Omg |\nabla \uu|^2,
 \label{StokesPE}
\end{equation}
where $\pcon=\frac23+\eps$.
\end{theorem}

This theorem allows one to see that \qref{NSE5} is {\em fully dissipative}.
To begin to see why, recall that the Laplace operator
$\Delta\colon H^2(\Omg)\cap H^1_0(\Omg)\to L^2(\Omg)$ is an
isomorphism, and note that $\grad\ps$ is determined by $\lap\uu$ via
\begin{equation}\label{ps-proj}
\nabla \ps = (I-\PP - \QQ)\lap\uu, \qquad \QQ:=\nabla \div \Delta^{-1}.
\end{equation}
Equation \qref{NSE5} can then be written
\begin{eqnarray}
 \pa_t \uu + \PP( \udotgrad \uu - \ff)
 &=& \nu (\PP+\QQ) \lap\uu  \nonumber\\
 &=& \nu \lap \uu  - \nu(I-\PP-\QQ)\lap\uu.  \label{NSE7}
\end{eqnarray}
Theorem~\ref{T.main} will allow us to regard the last term as a controlled
perturbation.

We can take $\lap\uu=\gg$ arbitrary in $L^2(\Omg,\R^N)$ and
reinterpret Theorem~\ref{T.main} as follows.
The last term in \qref{StokesPE}
can be interpreted as the squared norm of $\uu$ in $H^1_0(\Omg,\R^N)$,
giving the norm of $\gg$ in the dual space $H^{-1}(\Omg,\R^N)$.
Thus the conclusion of Theorem~\ref{T.main} is equivalent to the following
estimate, which says that $I-\PP$ is approximated
by the bounded operator $\QQ:L^2(\Omg,\R^N)\to\grad H^1(\Omg)$:
\begin{corollary}\label{C.1}
For all vector fields $\gg\in L^2(\Omg,\R^N)$ we have
\begin{equation} \label{E.Pc-Q}
\|(I-\PP - \QQ)\gg \|_{L^2}^2
\le
\pcon\|\gg\|_{L^2}^2 + C\|\gg\|_{H^{-1}}^2.
\end{equation}
\end{corollary}

There are several different ways to interpret the Stokes pressure
as we have defined it.  In this vein we make a few further
observations.  First, note that $\PP\nabla(\div \uu)=0$ for all
$\uu$ in $\uspace$, since $\div\uu$ lies in $H^1(\Omg)$.  
Then 
\begin{equation}\label{crossP}
\nabla \ps =(I-\PP)(\Delta\uu-\nabla(\div u)).
%\nabla \ps =-(I-\PP)(\nabla\times\nabla\times\uu).
\end{equation}
Now $\bau:=\Delta\uu-\nabla(\div\uu)$ has zero divergence
in the sense of distributions and is in
$L^2(\Omg,\R^N)$, so $\bau$ lies in the space $\Hdiv$ consisting of
vector fields in $L^2(\Omg,\R^N)$ with divergence in $L^2(\Omg)$.
By consequence, $\Delta\ps=0$ in the sense of distributions and so
$\nabla\ps$ is in $\Hdiv$ also. By a well-known trace theorem (see
\cite{GR}, theorem 2.5), the normal components of $\bau$ and
$\nabla\ps$ belong to the Sobolev space $H^{-1/2}(\Gamma)$, and
from the definition of $\PP$ we have
\begin{equation}\label{ps-wk1}
0 = \int_\Omg (\nabla\ps-\bau)\cdot\nabla\phi
=\int_\Gamma \phi\nn\cdot(\nabla\ps-\bau)
\end{equation}
for all $\phi\in H^1(\Omg)$.
So $\ps$ is determined as the zero-mean solution of the 
Neumann boundary-value problem
\begin{equation}\label{NeumannP}
\Delta\ps=0 \quad\mbox{in $\Omg$},\qquad
\nn\cdot\nabla\ps=\nn\cdot(\lap-\grad\div)\uu \quad\mbox{on $\Gamma$}.
\end{equation}

Furthermore, in two and three dimensions, we have
\begin{equation}\label{curlP}
\nabla \ps =-(I-\PP)(\nabla\times\nabla\times\uu)
\end{equation}
due to the identity
$\nabla\times\nabla\times\uu=-\Delta\uu+\nabla(\div u)$. 
Green's formula yields
\begin{equation}\label{mean0}
\int_\Gamma \nn\cdot(\grad\times\grad\times\uu)\phi
=\int_\Omg(\nabla\times\nabla\times\uu)\cdot\grad\phi
= -\int_\Gamma (\nabla\times\uu)\cdot(\nn\times\grad\phi)
\end{equation}
and so $\ps$ (with zero average) is determined 
through the weak formulation \cite{JL}
\begin{equation}\label{ps-wk}
\int_\Omg \grad\ps\cdot\grad\phi = \int_\Gamma
(\nabla\times\uu)\cdot(\nn\times\grad\phi)
\quad\text{for all $\phi\in H^1(\Omg)$.}
\end{equation}
(Note that $\nabla\times\uu\in H^{1/2}(\Gam,\R^N)$, and
$\nn\times\grad\phi\in H^{-1/2}(\Gam,\R^N)$ by a standard trace
theorem \cite[Theorem 2.11]{GR}, since $\grad\phi$ lies in
$\Hcurl$, the space of vector fields in $L^2(\Omg,\R^N)$ with curl in
$L^2$.) 

As indicated by \qref{NeumannP} or \qref{ps-wk}, 
the Stokes pressure is generated by
the {\it tangential part of vorticity at the boundary}.
In the whole space $\R^N$ or in the case of a periodic
box without boundary, the Helmholtz projection is exactly given via
Fourier transform by $\PP=I-\QQ$ and the Stokes pressure vanishes.
Essentially, the Stokes pressure supplies the correction to this formula
induced by the no-slip boundary conditions.
By consequence, the results of the present paper should have nothing to
do with the global regularity question for the three-dimensional
Navier-Stokes equations.  But as we have mentioned, many important
physical phenomena modeled by the Navier-Stokes equations involve
boundaries and boundary-layer effects, and it is exactly here where
the Stokes pressure should play a key role.

The unconstrained formulation \qref{NSE5} is not without antecedents
in the literature. Orszag et al.~\cite{OID} used the boundary condition
in \qref{NeumannP} as a way of enforcing consistency for a Neumann problem
in the context of the projection method.
After the results of this paper were completed, we found that
the formulation \qref{NSE5} is exactly equivalent to one 
studied by Grubb and Solonnikov \cite{GS1,GS2}.
These authors also study several other types of boundary conditions,
and argue that this formulation is parabolic in a nondegenerate sense.
They perform an analysis based on a general theory of parabolic
pseudo-differential initial-boundary value problems, 
and also show that for strong solutions,
the divergence satisfies a diffusion equation with Neumann boundary conditions.

Due to our Theorem~\ref{T.main}, we can treat the
Navier-Stokes equations in bounded domains simply
as a perturbation of the vector diffusion equation 
$\pa_t\uu = \nu\Delta\uu$, regarding
both the pressure and convection terms as dominated by the viscosity term.
This stands in contrast to the usual approach that regards
the Navier-Stokes equations as a perturbation of the Stokes system
$\pa_t\uu = \nu\Delta\uu-\grad p$, $\div\uu=0$.
Discussing this usual approach to analysis, Tartar \cite[p.\ 68]{Ta2} comments
\begin{quote}{``The difficulty comes from the fact that one does not have adequate
boundary conditions for $p$. $\ldots$[S]ending the nonlinear term to
play with $f$, one considers the Navier-Stokes equations as a
perturbation of Stokes equation, and this is obviously not a good
idea, but no one has really found how to do better yet.''}
\end{quote}
By way of seeking to do better, in this paper
we exploit Theorem~\ref{T.main} in a number of ways.  In particular, 
we develop a simple proof of local-in-time existence and uniqueness for strong
solutions of the unconstrained formulation \qref{NSE6} and
consequently for the original Navier-Stokes equations,
based upon demonstrating the {\em unconditional stability} of a simple
time-discretization scheme with explicit time-stepping for the
pressure and nonlinear convection terms and that is implicit only in the
viscosity term.

The discretization that we use is related to a class of extremely
efficient numerical methods for incompressible flow \cite{Ti96,Pe,JL,GuS}.  
Thanks to the explicit treatment of the convection and pressure terms,
the computation of the momentum equation is completely decoupled from
the computation of the kinematic pressure Poisson equation used to
enforce incompressibility.  No stationary Stokes solver is necessary
to handle implicitly differenced pressure terms.  
For three-dimensional flow in a general domain, the computation of
incompressible Navier-Stokes dynamics is basically reduced to solving
a heat equation and a Poisson equation at each time step.  
This class of methods is very flexible and can be used with all kinds of spatial
discretization methods \cite{JL}, including finite difference,
spectral, and finite element methods.  The stability properties we
establish here should be helpful in analyzing these methods.

Indeed, we will show below that our stability analysis 
easily adapts to proving unconditional
stability and convergence for corresponding fully discrete finite-element 
methods with $C^1$ elements for velocity and $C^0$ elements for pressure.
It is important to note that we impose {\em no} inf-sup compatibility
condition between the finite-element spaces for velocity and pressure.
The inf-sup condition (also known as the Ladyzhenskaya-Babu\v ska-Brezzi
condition) has long been a central foundation for finite-element
methods for all saddle-point problems including the stationary Stokes
equation. Its beautiful theory is a masterpiece documented in many
finite-element books. 
In the usual approach, the inf-sup condition serves
to force the approximate solution to stay close to the divergence-free
space where the Stokes operator $\PP \lap$ is dissipative.  
However, due to the fully dissipative nature of the unconstrained
formulation \qref{NSE6} which follows as a consequence of 
Theorem~\ref{T.main}, as far as our stability analysis in section~\ref{S.fem}  
is concerned, the finite-element spaces for velocity and pressure can be 
completely unrelated.

%The organization of this paper is as follows. 
The proof of Theorem~\ref{T.main}
will be carried out in section~\ref{S.pressure}. Important ingredients in the
proof are: (i) an estimate near the boundary that is 
related to boundedness of the
Neumann-to-Dirichlet map for boundary values of harmonic functions
--- this estimate is proved in section~\ref{S.N2D}, see Theorem~\ref{T.N2D}; 
and (ii) a representation formula for the
Stokes pressure in terms of a part of velocity near and parallel
to the boundary. In section~\ref{S.N2D} we also describe the space $\grad\SP$
of all possible Stokes pressure gradients (i.e., the range of $I-\PP-\QQ$).
In $\R^3$ it turns out that this is the space of
square-integrable vector fields that are {\em simultaneously
gradients and curls} (see Theorem~\ref{T.spgrad} in section~\ref{S.SP} below).  

In section~\ref{S.stab} we establish the unconditional stability of the
time-discretiza\-tion scheme, and in section~\ref{S.exist} we use this to study
existence and uniqueness for strong solutions with no-slip
boundary conditions. In section~\ref{S.fem} we adapt the stability analysis to
prove the unconditional stability and convergence of corresponding 
$C^1/C^0$ finite-element methods.

In section~\ref{S.semig} we show that Theorem~\ref{T.main} also allows one
to treat the linearized equations (an unconstrained version of the
Stokes system) easily by analytic semigroup theory.  We deal with
non-homogeneous boundary conditions in section~\ref{S.nonhom}.  From
these results, in section~\ref{S.stokes} we deduce an apparently new
result for the linear Stokes system, namely an isomorphism theorem
between the solution space and a space of data for non-homogeneous
side conditions in which only the average flux through the boundary
vanishes.

\section{Integrated Neumann-to-Dirichlet estimates in tubes}
\label{S.N2D}

\subsection{Notation}
Let $\Omg\subset\R^N$ be a bounded domain with $C^3$ boundary
$\Gam$. For any $\xx\in\Omg$ we let $\Phi(\xx) = \dist(x,\Gamma)$
denote the distance from $x$ to $\Gam$. For any $s>0$ we denote
the set of points in $\Omg$ within distance $s$ from $\Gamma$ by
\begin{equation} \label{E.Omgs}
\Omg_s = \{ \xx \in \Omg \mid \Phi(\xx) < s \},
\end{equation}
and set $\Omg_s^c = \Omg \backslash \Omg_s$ and $\Gam_s = \{ \xx
\in \Omg \mid \Phi(\xx) = s \}$. Since $\Gamma$ is $C^3$ and
compact, there exists $s_0>0$ such that $\Phi$ is $C^3$ in
$\Omg_{s_0}$ and its gradient is a unit vector, with $|\grad
\Phi(\xx)| = 1$ for every $\xx\in\Omg_{s_0}$. We let
\begin{equation} \label{E.nn}
   \nn(\xx) = -\grad \Phi(\xx),
\end{equation}
then $\nn(\xx)$ is the outward unit normal to $\Gam_s=\pa\Omg_s^c$
for $s=\Phi(\xx)$, and $\nn \in C^2(\bOmg_{s_0},\R^N)$.

We let $\<f,g\>_\Omg = \int_{\Omg}fg$ denote the $L^2$ inner product
of functions $f$ and $g$ in $\Omg$, and let $\| \cdot \|_\Omg$ denote
the corresponding norm in $L^2(\Omg)$.  We drop the subscript on the
inner product and norm when the domain of integration is understood in
context.

\subsection{Statement of results}

Our strategy for proving Theorem~\ref{T.main} crucially involves an
integrated Neumann-to-Dirichlet--type estimate for harmonic functions
in the tubular domains $\Omg_s$ for small $s>0$.

The theorem below contains two estimates of this type. 
The first, \qref{plap1} in part (i), can be obtained from
a standard Neumann-to-Dirichlet estimate for functions harmonic
in $\Omg$, of the form 
\begin{equation}\label{N2D1}
\beta_0 \int_{\Gam_r} |(\tanproj)\grad p|^2 \le
\int_{\Gam_r} |\nder p|^2 ,
\end{equation}
by integrating over $r\in(0,s)$, provided one shows that $\beta_0>0$
can be chosen independent of $r$ for small $r>0$. 
On the first reading, the reader is encouraged to take \qref{plap1} for
granted and proceed directly to section 3.2 at this point; it is only
necessary to replace \qref{upest1} in the proof of
Theorem~\ref{T.main} by the corresponding result from \qref{plap1} to
establish that the estimate in Theorem~\ref{T.main} is valid for {\em
some} $\pcon<1$ depending upon $\Omg$.

The second estimate, in part (ii), will be used with $\beta_1$ close
to $1$ to establish the full
result in Theorem~\ref{T.main} for any number $\pcon$ greater than
$\frac23$, independent of the domain.

\begin{theorem}\label{T.N2D}
Let $\Omg$ be a bounded domain with $C^3$ boundary.
(i) There exists $\beta_0>0$ such that for sufficiently small $s>0$,
whenever $p$ is a harmonic function in $\Omg_s$ we have
\begin{equation}\label{plap1}
\beta_0 \int_{\Omg_s} |(\tanproj)\grad p|^2 \le 
\int_{\Omg_s} |\nder p|^2.
\end{equation}
(ii) Let $\beta_1<1$. Then for any sufficiently small $s>0$,
whenever $p$ is a harmonic function in $\Omg_s$ and 
$p_0$ is constant on each component of $\Omg_s$, we have
\begin{equation}
\beta_1\int_{\Omg_s} |(\tanproj)\grad p|^2 \le 
\int_{\Omg_s} |\nder p|^2 + \frac{24}{s^2}\int_{\Omg_s} |p-p_0|^2.
\end{equation}
\end{theorem}

Our proof is motivated by the case of slab domains with
periodic boundary conditions in the transverse directions. In this case
the analysis reduces to estimates for Fourier series expansions
in the transverse variables.  For general domains, 
the idea is to approximate $-\lap$ in thin tubular domains $\Omg_s$ 
by the Laplace-Beltrami operator on $\Gam\times(0,s)$.  
This operator has a direct-sum structure, and 
we obtain the integrated Neumann-to-Dirichlet--type estimate by
separating variables and expanding in series of eigenfunctions of the
Laplace-Beltrami operator on $\Gamma$.
For basic background in Riemannian geometry and the Laplace-Beltrami
operator we refer to \cite{Au} and \cite{Ta}.

\subsection{Harmonic functions on $\Gam\times(0,s)$}
{\bf Geometric preliminaries.}
We consider the manifold $\G=\Gam\times\I$ with $\I=(0,s)$
as a Riemannian submanifold of $\R^N\times\R$ with boundary
$\pa\G=\Gam\times\{0,s\}$.
We let $\gamma$ denote the metric on $\Gamma$ induced from $\R^N$, 
let $\iota$ denote the standard Euclidean metric on $\I$,
and let $g$ denote the metric on the product space $\G$.
Any vector $\ba$ tangent to $\G$ at $z=(y,r)$ 
has components $\ba_\Gam$ tangent to $\Gam$ at $y$ and $\ba_\I$ tangent to
$\I$ at $r$. For any two such vectors $\ba$ and $\bb$, we have
\begin{equation}
g(\ba,\bb) = 
\gamma(\ba_\Gam,\bb_\Gam)+\iota(\ba_\I,\bb_\I).
\label{g-GI1}
\end{equation}

Given a $C^1$ function $z=(y,r)\mapsto f(y,r)$ on $\G$, its gradient
$\gradg f$ at $z$ is a tangent vector to $\G$ determined from 
the differential via the metric, through requiring
\begin{equation}
g(\gradg f,\ba) = df\cdot\ba \quad\mbox{for all $\ba\in T_z\G$}.
\label{g-GI2}
\end{equation}
By keeping $r$ fixed, the function $y\mapsto f(y,r)$ determines 
the gradient vector $\gradG f$ tangent to $\Gam$ in similar fashion,
and by keeping $y$ fixed, the function $r\mapsto f(y,r)$ determines
the gradient vector $\gradI f$ tangent to $\I$.
These gradients are also the components of $\gradg f$:
\[
(\gradg f)_\Gam = \gradG f, \qquad (\gradg f)_\I = \gradI f.
\]

If $u=(u^1,\ldots,u^{N-1})\mapsto y=(y^1,\ldots,y^N)$ 
is a local coordinate chart for $\Gam$, 
the metric is given by $\gam_{ij}\,du^i\,du^j$ 
(summation over repeated indices implied) with matrix elements 
\[
\gam_{ij} = 
\frac{\pa y^k}{\pa u^i}
\frac{\pa y^k}{\pa u^j} .
\]
For $\I\subset\R$ the identity map serves as coordinate chart.
In these coordinates the tangent vectors are written 
(in a form that aids in tracking coordinate changes) as
\begin{equation}
\gradG f =  \gam^{ij}\frac{\pa f}{\pa u^i}\frac{\pa}{\pa u^j},
\qquad
\gradI f = \frac{\pa f}{\pa r}\frac{\pa}{\pa r}.
\end{equation}
As usual, the matrix $(\gam^{ij})=(\gam_{ij})^{-1}$.
Given two $C^1$ functions $f, \tf$ on $\G$, 
\begin{equation}
\gam(\gradG f,\gradG\tf) = \gam^{ij}
\frac{\pa f}{\pa u^i}
\frac{\pa\tf}{\pa u^j},
\qquad
\iota(\gradI f,\gradI\tf) = \frac{\pa f}{\pa r}\frac{\pa\tf}{\pa r}.
\label{gi1}
\end{equation}

In these coordinates, the (positive) Laplace-Beltrami operators on
$\Gam$ and $\I$ respectively take the form 
\begin{equation}\label{LBdef1}
\lapG f = -\frac{1}{\sqrt{\gam}} \frac{\pa}{\pa u_i}
\left(\sqrt{\gam}\gam^{ij}\frac{\pa}{\pa u_j}f\right) ,
\qquad
\lapI f = -\frac{\pa^2}{\pa r^2} f ,
\end{equation}
where $\sqrt{\gam} = \sqrt{\det(\gam_{ij})}$ is the
change-of-variables factor for integration on $\Gam$ --- 
if a function $f$ on $\Gam$ is
supported in the range of the local coordinate chart then
\begin{equation}
\int_\Gam f(y)\,dS(y) = \int_{\R^{N-1}} f(y(u)) \,\sqrt{\gam}\,du.
\end{equation}
(Since orthogonal changes of coordinates in $\R^N$ and $\R^{N-1}$
leave the integral invariant, one can understand $\sqrt{\gam}$ as the
product of the singular values of the matrix $\pa y/\pa u$.)

Whenever $f\in H^1(\Gam)$ and $\tf\in H^2(\Gam)$,
one has the integration-by-parts formula
\begin{equation}
 \int_\Gam {f\lapG\tf} =
 \int_\Gam \gam(\gradG f,\gradG\tf)
.% = \int_0^s\int_\Gam (\lapg f)\tf \,dS(y)\,dr.
\label{intparts}
\end{equation}
One may extend $\lapG$ to be a map from $H^1(\Gam)\to H^{-1}(\Gam)$ by
using this equation as a definition of $\lapG \tf$ as a functional on
$H^1(\Gam)$. In standard fashion \cite{Ta}, one finds that 
$I+\lapG:H^1(\Gam)\to H^{-1}(\Gam)$ is an isomorphism, and 
that $(I+\lapG)^{-1}$ is a compact self-adjoint operator on
$L^2(\Gam)$, hence $L^2(\Gam)$ admits an orthonormal basis of 
eigenfunctions of $\lapG$.  Since the coefficient functions in
\qref{LBdef1} are $C^1$, 
standard interior elliptic regularity results (\cite[Theorem 8.8]{GT},
\cite[p. 306, Proposition 1.6]{Ta})
imply that the eigenfunctions belong to $H^2(\Gam)$.
We denote the eigenvalues of $\lapG$ by $\nu_k^2$, $k=1,2,\ldots$, with
$0=\nu_1\le\nu_2\le\ldots$ where $\nu_k\to\infty$ as $k\to\infty$, 
and let $\psi_k$ be corresponding eigenfunctions 
forming an orthonormal basis of $L^2(\Gam)$.  
If $\lapG \psi=0$ then $\psi$ is
constant on each component of $\Gam$, so if $m$ is the number of
components of $\Gam$, then $0=\nu_m<\nu_{m+1}$.

In the coordinates $\hat u=(u,r)\mapsto z=(y,r)$ for $\G$, the metric 
$g$ takes the form $\gam_{ij}\,du^i\,du^j+ dr^2$, and 
the Laplace-Beltrami operator $\lapg=\lapG+\lapI$.
Similar considerations as above apply to 
$\lapg$, except $\G$ has boundary. 
Whenever $f\in H^1_0(\G)$ and $\tf\in H^2(\G)$
we have
\begin{equation}
 \int_\G {f\lapg\tf} =
 \int_\G g(\gradg f,\gradg\tf).
% = \int_0^s\int_\Gam (\lapg f)\tf \,dS(y)\,dr.
\label{intparts2}
\end{equation}
One extends $\lapg$ to map $H^1(\G)$ to $H^{-1}(\G)$ by
using this equation as a definition of $\lapg \tf$ as a functional on
$H^1_0(\G)$.

We introduce notation for $L^2$ inner products and norms on $\G$ as
follows:
\begin{align}
\ip{f,\tf}_\G &= \int_\G f\tf 
&\|f\|_\G^2 &= \int_\G |f|^2,
\\ \label{ipg1}
\ip{\gradG f,\gradG \tf}_\G &= \int_\G \gam(\gradG f,\gradG\tf) ,
&\|\gradG f\|_\G^2 &= \int_\G \gam(\gradG f,\gradG f),
\\ \label{ipg2}
\ip{\gradI f,\gradI \tf}_\G &= \int_\G (\pa_rf)(\pa_r\tf) ,
&\|\gradI f\|_\G^2 &= \int_\G (\pa_rf)^2 ,
\end{align}
\begin{align} \label{ipg3}
\ip{\gradg f,\gradg\tf}_\G &= 
\int_\G g(\gradg f,\gradg \tf) = 
\ip{\gradG f,\gradG\tf}_\G +\ip{\gradI f,\gradI\tf}_\G,
\\
\|\gradg f\|_\G^2 &= 
\int_\G g(\gradg f,\gradg f) = 
\|\gradG f\|_\G^2 + \|\gradI f\|_\G^2 .
\label{Ginp}
\end{align}

\begin{lemma} \label{L.N2D}
Suppose $f\in H^1(\G)$ and $\lapg f=0$ on $\G=\Gam\times\I$ where $\I=(0,s)$. 
Then, (i) there exists
$\hat\beta_0\in(0,1)$ independent of $f$ such that
\begin{equation}
\hat\beta_0 \|\gradG f\|_\G^2 \le \|\gradI f\|_\G^2 ,
\end{equation}
and (ii) 
\begin{equation}
\|\gradG f\|_\G^2 \le \|\gradI f\|_\G^2 + \frac{12}{s^2}
\|f-f_0\|_\G^2
\end{equation}
whenever $f_0$ is constant on $\Gam_i\times(0,s)$ for every component
$\Gamma_i$ of $\Gam$.
\end{lemma}

\noindent{\bf Proof:}
Suppose $\lapg f=0$ on $\G$. 
Since the coefficient functions in
\qref{LBdef1} are $C^1$, the aforementioned interior elliptic regularity results
imply that that $f\in \Htwoloc(\G)$. 
For any $r\in(0,s)$, fixing $r$ yields a trace of $f$ in
$H^1(\Gam)$, and as a function of $r$, 
we can regard $f=f(y,r)$
as in the space $L^2([a,b],H^2(\Gam))\cap H^2([a,b],L^2(\Gam))$
for any closed interval $[a,b]\subset(0,s)$.
Now, for each $r$ we have the $L^2(\Gam)$-convergent expansion
\begin{equation}\label{fsum}
f(y,r)=\sum_k \hat f(k,r) \psi_k(y)
\end{equation}
where
\begin{equation}\label{hatf}
\hat f(k,r) = \int_{\Gam} f(y,r)\psi_k(y) \,dS(y) .
%= \int_{\Gam_r} f(x)\psi_k(x+r\nn(x))\sqrt{g}\,dS(x) .
\end{equation}
For each $k\in\N$, the map $r\mapsto \hat f(k,r)$ is in $\Htwoloc(0,s)$
and 
\begin{equation}\label{E.hatfr}
\pa_r \hat f(k,r) 
%= \pa_r \int_\Gam f(y,r)\psi_k(y)
= \int_\Gam \pa_r f(y,r)\psi_k(y) \,dS(y) .
%\ip{f_r,\psi_k}_\Gam = \ip{(\pa_rf)_r,\psi_k}_\Gam.
\end{equation}
For any smooth $\xi\in C^\infty_0(0,s)$, taking
$\tf(y,r)=\psi_k(y)\xi(r)$ we compute that
\begin{equation}\label{gradtf}
\gradG\tf = \xi(r) \gradG\psi_k,\qquad \pa_r\tf= \psi_k \pa_r\xi,
\end{equation}
and so by \qref{intparts2}, \qref{g-GI1}, and \qref{intparts}, we have
\begin{align}
%0&= \ip{\lapg f,\tf}_\G = 
0&= \int_\G (\lapg f)\tf = 
\int_\I\int_\Gam \left(
\gam(\gradG f,\gradG\tf)+(\pa_rf)(\pa_r\tf)\right)
\nonumber\\& = 
\int_\I \xi(r) \int_\Gam 
\gam(\gradG f,\gradG\psi_k)
+ \int_\I(\pa_r\xi)\int_\Gam (\pa_rf)\psi_k
\nonumber\\& =
\int_\I \xi(r) \int_\Gam f \lapG\psi_k
+ \int_\I (\pa_r\xi) \pa_r\hat f(k,r)
\nonumber\\&=
\int_0^s \left(\xi(r) \nu_k^2 \hat f(k,r)
+ (\pa_r\xi)\pa_r\hat f(k,r)\right)\,dr .
\label{weak0}
\end{align}
Therefore $\hat f(k,\cdot)$ is a weak solution of 
$\pa_r^2\hat f=\nu_k^2\hat f$ in $\Htwoloc(0,s)$ and hence is $C^2$ and
it follows that whenever $\nu_k\ne0$, there exist $a_k$, $b_k$ such
that 
\begin{equation}
\hat f(k,r) = a_k\sinh \nu_k\tau + b_k\cosh \nu_k\tau, \qquad
\tau=r-s/2.
\end{equation}

Now
\begin{align}
\|f\|_\G^2 &= \sum_k \int_0^s |\hat f(k,r)|^2 \,dr,
\\
\|\gradG f\|_\G^2 &= \sum_k \int_0^s |\nu_k \hat f(k,r)|^2 \,dr,
\\
\|\gradI f\|_\G^2 &= \sum_k \int_0^s |\pa_r \hat f(k,r)|^2 \,dr.
\end{align}
Let $\gamma_k=\int_{-s/2}^{s/2} \sinh^2\nu_k \tau\,d\tau$.
Then $\gamma_k$ increases with $k$, and 
\begin{equation}
\gamma_k+s = \int_{-s/2}^{s/2} \cosh^2\nu_k \tau\,d\tau
\ge \int_{-s/2}^{s/2} (1+\nu_k^2\tau^2)\,d\tau \ge \frac{\nu_k^2 s^3}{12}.
\end{equation}
Whenever $\nu_k\ne0$ we get
\begin{align}
\int_0^s |\hat f(k,r)|^2\,dr &= 
 |a_k|^2\gamma_k+|b_k|^2 (\gamma_k+s),
\\
\int_0^s |\pa_r \hat f(k,r)|^2\,dr &= 
\nu_k^2( |a_k|^2(\gamma_k+s)+|b_k|^2 \gamma_k),
\end{align}
and since $\hat\beta_0(\gamma_k+s)\le \gamma_k$ where
$\hat\beta_0=\gamma_{m+1}/(\gamma_{m+1}+s)$, it follows
\begin{align}
\hat\beta_0\int_0^s |\nu_k \hat f(k,r)|^2 \,dr &\le 
\int_0^s |\pa_r \hat f(k,r)|^2 \,dr,
\\
\int_0^s |\nu_k \hat f(k,r)|^2 \,dr &\le \int_0^s |\pa_r \hat f(k,r)|^2 \,dr
+\frac{12}{s^2} \int_0^s |\hat f(k,r)|^2\,dr.
\end{align}
The results in (i) and (ii) follow by summing over $k$. 
$\square$

\subsection{Global coordinates on $\Gam\times(0,s)$}

It will be important for comparison with the Laplacian on $\Omg_s$
to coordinatize $\G$ for small $s>0$ {\em globally} 
via the coordinate chart $\Omg_s\to \G$ given by
\begin{equation}\label{g-chart}
x\mapsto z = (y,r) = (x+\Phi(x)\nn(x),\Phi(x)).
\end{equation}
In these coordinates, the metric on $\G$ that is inherited from
$\R^{N+1}$ has the representation $g_{ij}\,dx^i\,dx^j$ with matrix
elements given by 
\begin{equation}
g_{ij} = 
\frac{\pa z^k}{\pa x^i}
\frac{\pa z^k}{\pa x^j} =
\frac{\pa y^k}{\pa x^i}
\frac{\pa y^k}{\pa x^j} +
\frac{\pa r}{\pa x^i}
\frac{\pa r}{\pa x^j} .
\end{equation}
Let us write $\pa_i=\pa/\pa x^i$ and let 
$\nabla f=(\pa_1 f,\ldots,\pa_N f)$ 
denote the usual gradient vector in $\R^N$. The components of $\nn$
are $n_i = -\pa_i\Phi$ and so $\pa_in_j=\pa_jn_i$, meaning the matrix
$\grad\nn$ is symmetric. Since $|\nn|^2=1$
we have $n_i\pa_jn_i=0=n_i\pa_in_j$.  Then the $N\times N$ matrix
\begin{equation}
\frac{\pa y}{\pa x} = \tanproj+\Phi\grad\nn = 
(\tanproj)(I+\Phi\grad\nn)(\tanproj) ,
\end{equation}
and the matrix
\begin{equation}\label{E.met2}
G = (g_{ij}) = (\tanproj)(I+\Phifac)^2(\tanproj)+\nproj=
(I+\Phifac)^2 .
\end{equation}

With $\sqrt{g}=\sqrt{\det G}$, 
the integral of a function $f$ on $\G$ in terms of these coordinates 
is given by 
\begin{equation}
\int_\G f = \int_\Omgs f\,\sqrt{g}\,dx .
\label{intG}
\end{equation}
Given two $C^1$ functions $f$, $\tf$ on $\G$, 
we claim that the following formulae are valid in the coordinates
from \qref{g-chart}:
\begin{align}
g(\gradg f,\gradg \tf) &= (\grad f)^t G^{-1} (\grad\tf) = 
g^{ij}\pa_if\pa_j\tf ,
\label{eq-g}
\\
\gamma(\gradG f,\gradG \tf) &= 
(\nabla f)^t(\tanproj) G^{-1}(\tanproj) (\nabla\tf) ,
\label{eq-gG}
\\
\iota(\gradI f,\gradI \tf) &= (\nder f)(\nder\tf) = 
(\grad f)^t\nproj(\grad\tf) .
\label{eq-gI}
\end{align}
Of course \qref{eq-g} simply expresses the metric in the $x$-coordinates
from \qref{g-chart}.
To prove \qref{eq-gI}, first note that along any curve $\tau\mapsto
x(\tau)$ satisfying $\pa_\tau x = \nn(x)$ we find
$\pa_\tau\nn(x)=n_j\pa_jn_i=0$, so $\nn(x)$ is constant and the curve
is a straight line segment.
Hence in the chart from \qref{g-chart}, $\nn(x)=\nn(y)$ and we have
$x= y-r\nn(y)$. Given a $C^1$ function $f$ 
then, we find that in these $\Omgs$-coordinates,
\begin{equation}
\pa_r f(y,r) = (\pa_r x_j)(\pa_j f) =  n_j \pa_j f = \nder f,
\end{equation}
and \qref{eq-gI} follows from \qref{gi1}.
Finally, \qref{eq-gG} follows directly from \qref{eq-g} and
\qref{eq-gI} using \qref{g-GI1} --- since
$\nproj\grad\nn=0$ we have $\nproj G=\nproj$ so $\nproj=\nproj G^{-1}$
and hence
\begin{equation}
(\tanproj)G^{-1}(\tanproj) = G^{-1}-\nproj.
\end{equation}

\subsection{Proof of Theorem~\ref{T.N2D}}
Let $\beta_1<1$.  Suppose $\lap p=0$ in $\Omg_s$.
We may assume $p\in H^1(\Omgs)$ without loss of generality
by establishing the result in subdomains where $\Phi(x)\in(a,b)$ with
$[a,b]\subset(0,s)$ and taking $a\to0$, $b\to s$.
We write 
\[
p = p_1+p_2,
\]
where $p_1\in H^1_0(\Omg_s)$ is found by solving
a weak form of $\lapg p_1=\lapg p$:
\begin{equation}\label{p1eqn}
\ip{\gradg p_1,\gradg\phi}_\G = \ip{\gradg p,\gradg\phi}_\G
\quad\mbox{for all $\phi\in H^1_0(\Omg_s)$.}
\end{equation}
For small $s>0$, $G=(g_{ij})=I+O(s)$ and
$\sqrt{g}=1+O(s)$. Since $\ip{\grad p,\grad p_1}=0$, 
taking $\phi=p_1$ we have
\begin{equation}
\|\gradg p_1\|_\G^2 =
\int_{\Omg_s} (\grad p)^t (G^{-1}\sqrt{g} -I) \grad p_1\,dx
\le Cs \|\grad p\|_\Omgs \|\gradg p_1\|_\G ,
\end{equation}
where $C$ is a constant independent of $s$.
By Poincar\'e's inequality we also have 
\begin{equation}
\|p_1\|_\G^2 \le \frac{s^2}{\pi^2}\|\gradg p_1\|_\G^2
\end{equation}
since the eigenvalues of $\lapg$ on the product space 
$\Gam\times[0,s]$ with Dirichlet
boundary conditions all have the form $\mu=\nu_k^2+j^2\pi^2/s^2$
for $j,k\in\N$, so that $\mu\ge \pi^2/s^2$.

Let us first prove part (ii).  For $0<\eps<1$, using
\qref{eq-gG}, \qref{intG} and \qref{ipg1} we deduce
\begin{eqnarray}
\|(\tanproj)\grad p\|_\Omgs^2 &\le& (1+Cs)\|\gradG p\|_\G^2
\nonumber\\&\le& 
(1+Cs)\left( (1+\eps)
\|\gradG p_2\|_\G^2 + (1+\eps^{-1})\|\gradG p_1\|_\G^2 \right)
\nonumber\\&\le& 
(1+Cs)(1+\eps)\left( 
\|\gradG p_2\|_\G^2 + \eps^{-1}C^2s^2 \|\grad p\|_\Omgs^2 \right).
\label{p2est0}
\end{eqnarray}
Now $p_2=p -p_1$ satisfies $\lapg p_2=0$ in $\Omg_s$ and 
$p_2\in H^1(\G)$, hence
for any $p_0$ constant on each component of $\Omg_s$ we have
\begin{eqnarray}
\|\gradG p_2\|_\G^2 &\le& \|\gradI p_2\|_\G^2 +
\frac{12}{s^2}\|p_2-p_0\|_\G^2 ,  
\label{p2est1}
\\
\|\gradI p_2\|_\G^2 
&\le& (1+\eps)\|\gradI p\|_\G^2 + (1+\eps^{-1})\|\gradI p_1\|_\G^2 
\nonumber\\
&\le& (1+\eps)(1+Cs)\left(
\|\nder p\|_\Omgs^2 + \eps^{-1}C^2s^2\|\grad p\|_\Omgs^2 \right) ,
\label{p2est2}
\\
\frac{12}{s^2} \|p_2-p_0\|_\G^2 &\le& 
\frac{24}{s^2} \left( \|p -p_0\|_\G^2 + \|p_1\|_\G^2\right)
\nonumber\\
&\le& \frac{24}{s^2} \|p-p_0\|_\G^2+ \frac{24}{\pi^2} \|\gradg p_1\|_\G^2
\nonumber\\
&\le& \frac{24}{s^2}(1+Cs) \|p-p_0\|_\Omgs^2+ C^2s^2 \|\grad p\|_\Omgs^2
\label{p2est3}
\end{eqnarray} 
Presuming $Cs<\frac13\eps$,
assembling these estimates yields
\begin{eqnarray}
&& \|(\tanproj)\grad p\|_\Omgs^2 
\le (1+\eps)^4 \left(
\|\nder p\|_\Omgs^2 + \frac{24}{s^2}\|p-p_0\|_\Omgs^2 + \eps %\frac{\eps}3
\|\grad p\|_\Omgs^2\right) \nonumber\\
&&\quad \le\ (1+\eps)^5 \left(
\|\nder p\|_\Omgs^2 + \frac{24}{s^2}\|p-p_0\|_\Omgs^2 + \eps
\|(\tanproj)\grad p\|_\Omgs^2\right),
\end{eqnarray}
since
$|\grad p|^2= |\nder p|^2+|(\tanproj)\grad p|^2$. 
Fixing $\eps>0$ small so that $(1+\eps)^{-5}-\eps>\beta_1$ proves part (ii).

To prove part (i), instead of \qref{p2est1} we use 
\begin{equation}
\hat\beta_0\|\gradG p_2\|_\G^2 \le \|\gradI p_2\|_\G^2 
\end{equation}
(from part (i) of Lemma~\ref{L.N2D})
together with \qref{p2est0} and \qref{p2est2} and obtain 
\begin{eqnarray}
 \hat\beta_0 \|(\tanproj)\grad p\|_\Omgs^2 
&\le& (1+\eps)^4 \left(
\|\nder p\|_\Omgs^2 + \frac{2\eps}9
\|\grad p\|_\Omgs^2\right) \nonumber\\
&\le& (1+\eps)^5 \left(
\|\nder p\|_\Omgs^2 + \eps
\|(\tanproj)\grad p\|_\Omgs^2\right) .
\end{eqnarray}
Now taking $\eps>0$ so small that $\eps(1+\eps)^5<\hat\beta_0$ 
finishes the proof.
$\square$

\section{Analysis of the Stokes pressure}
\label{S.pressure}

The main purpose of this section is to prove Theorem~\ref{T.main}.  We
also describe the range of the map $\uu\mapsto\grad\ps$ from velocity
fields to Stokes pressure gradients.  For motivation for the proof of
Theorem~\ref{T.main}, the reader can proceed directly to
section~\ref{S.upeupa} at this point. Here, we first establish some
key preliminary results.

\subsection{An $L^2$ estimate}
The following $L^2$ estimate on the Stokes pressure will be used
to obtain the full result of Theorem~\ref{T.main} for arbitrary
$\pcon>\frac23$.  It is not needed to prove the weaker statement that
\qref{StokesPE} holds for some $\pcon<1$.

\begin{lemma} \label{L.pL2}
Let $\Omg\subset\R^N$ ($N\ge2$)
be any bounded connected domain with $C^{1,1}$ boundary. 
For any $\eps >0$, there is a constant $C\ge0$ so that for any 
$\uu \in \uspace$, the
associated Stokes pressure $\ps$ defined by $\qref{StokesP}$
with zero average satisfies
\begin{equation} \label{E.pL2.est}
  \|\ps \| \leq \eps \|\lap \uu \|
 + C \| \uu\|.
\end{equation}
%where $\bps = \frac{1}{|\Omg|} \int_{\Omg} \ps$.
\end{lemma}
{\bf Proof:} For any $\phi \in L^2(\Omg)$, define $\psi$ by
\begin{equation} \label{E.temp30}
  \lap \psi = \phi - \bar{\phi} ,\qquad \nn \cdot \grad \psi
  \big|_{\Gam} =0,
\end{equation}
where $\bar{\phi}$ is the average value of $\phi$ over $\Omg$.
Recall $\bps=0$. Then,
\begin{equation} \label{added-1}
  \<\ps , \phi \>  = \<\ps, \phi - \bar{\phi}\> = \< \ps, \lap
  \psi\>
   = -\< \grad \ps, \grad \psi \>
   % = -\< \curl \uu, \nn \times \grad \psi \>_{\Gam}
   .
\end{equation}
From \qref{ps-wk}, we know $ \< \grad \ps, \grad \psi \> = \<
\curl \uu, \nn \times \grad \psi \>_{\Gam}$ when $N=2$ or $3$. For
general $N$, 
using the notation $\pa_i:=\pa/\pa x_i$ and automatic summation
upon repeated indices, 
from \qref{NeumannP} we derive
\begin{equation} \label{added-2}
 \< \grad \ps, \grad \psi \> = \frac12 \int_\Gam (\pa_j u_i - \pa_i
 u_j)(n_j \pa_i \psi - n_i \pa_j \psi ).
\end{equation}
Plug \qref{added-2} into \qref{added-1}, take the absolute value
and use the trace theorem to get
%control the boundary integral coming from get
%$\big| \int_\Gam f g \big| \leq \|f\|_{L^2(\Gam)}
%\|g\|_{L^2(\Gam)}$. We get
\begin{equation}
  |\<\ps , \phi \>|
\leq c_0\|\grad\uu\|_{L^2(\Gam)} \|\grad\psi\|_{L^2(\Gam)}
%   \leq
%   \frac12 \| \pa_j u_i - \pa_i u_j \|_{L^2(\Gam)} \| n_j \pa_i \phi -
%   n_i \pa_j \phi \|_{L^2(\Gam)}
   \leq c_1 \|\uu\|_{H^{3/2}(\Omg)} \| \psi \|_{H^{3/2}(\Omg)}.
   \label{E.temp31}
\end{equation}
By the regularity theory for Poisson's equation \qref{E.temp30},
\begin{equation} \label{E.temp33}
\|\psi\|_{H^{3/2}(\Omg)} \leq c_2 \|\phi - \bar{\phi}\| \leq
c_2\|\phi\|.
\end{equation}
By a standard interpolation theorem, for any $\delta >0$, there is
a constant $c$, so 
\begin{equation} \label{E.temp32}
\|\uu\|_{H^{3/2}(\Omg)} \leq \delta \|\lap \uu\| + c \| \uu\| .
\end{equation}
Plugging \qref{E.temp33} and \qref{E.temp32} into \qref{E.temp31}, we
get
\begin{equation} \label{E.temp34}
|\<\ps , \phi \>| \leq \( \delta \|\lap \uu\| + c \| \uu\|\)
c_1c_2 \|\phi\|.
\end{equation}
Thus,
\begin{equation}
  \|\ps \| = \sup_{\phi \in L^2}\frac{|\<\ps , \phi
  \>|}{\| \phi\|} \leq \eps\|\lap \uu\| + cc_1c_2\|\uu\|.
  \qquad \square
\end{equation}

\subsection{Identities at the boundary}

A key part of the proof of Theorem~\ref{T.main} involves boundary
values of two quantities that involve the decomposition of
$\uu=(\tanproj)\uu+\nproj\uu$ into parts parallel and normal to the boundary.
Our goal in this subsection is to prove the following.
\begin{lemma} \label{L.uH2}
Let $\Omg \subset \R^N$ be a bounded domain with boundary
$\Gam$ of class $C^3$. Then for any $\uu \in H^2(\Omg, \R^N)$ with
$\uu|_\Gam =0 $, the following is valid on $\Gam$:
 \begin{description}
   \item (i) $\div \( (\tanproj)\uu \)=0 $ in $H^{1/2}(\Gam)$.
   \item (ii) $ \nn \cdot  ( \lap - \grad \div) \( \nproj \uu \) = 0 $
   in $H^{-1/2}(\Gam)$.
 \end{description}
\end{lemma}
The proof will reduce to the case $\uu\in
C^2(\bOmg,\R^N)$, due to the following density result.
\begin{lemma} \label{L.C2density}
Let $\Omg \subset \R^N$ be a bounded domain with boundary
$\Gam$ of class $C^{2,\alpha}$ where $0<\alpha<1$.
 Then for any $\uu \in \uspace$, there exists a
 sequence $\uu_k \in C^{2,\alpha}(\bOmg)$ such that $\uu_k|_\Gam =0$
 and $\| \uu_k -\uu \|_{H^2(\Omg)} \rightarrow 0$.
\end{lemma}
{\bf Proof:} Define $\ff = \lap\uu.$ Since $\ff \in L^2$, we can
find a sequence $\ff_k \in C^1(\bOmg)$ so that $\|\ff_k
-\ff\|_{L^2} \rightarrow 0.$ Construct $\uu_k$ by solving
$$
\lap \uu_k =\ff_k, \qquad \uu_k \big|_{\pa \Omg} = 0.
$$
Classical elliptic regularity theory in H\"older spaces (see \cite{GT},
theorem 15.13) says that a unique $\uu_k$ exists and is in
$C^{2,\alpha}(\bOmg)$. By standard regularity theory
in Sobolev spaces,
\[
\| \uu_k - \uu \|_{H^2} \leq C\|\ff_k -\ff\|_{L^2} \rightarrow 0.
\qquad\square
\]

\noind {\bf Proof of Lemma~\ref{L.uH2}:}
To begin, recall $\nn = -\grad \Phi$.
Equality of mixed partial derivatives yields
%\begin{equation} \label{E.n.id2}
$\pa_j n_i= \pa_i n_j$
%\end{equation}
for all $i, j=1,\ldots,N$. Together with the fact $n_in_i=1$, we infer
that for small $s>0$, throughout $\Omg_{s}$ we have
\begin{equation} \label{E.n.id1}
  n_i \pa_j n_i =0 \quad\text{and}\quad  n_i \pa_i n_j=0.
\end{equation}
(i) First, for any $f\in C^1(\bOmg)$, if $f=0$ on $\Gam$
then $\nabla f\parallel\nn$ on $\Gam$, which means
  \begin{equation} \label{E.n.id4.0}
  (\tanproj)\grad f =0, \quad\mbox{or}\quad
    (\pa_k - n_kn_j\pa_j)f=0 \quad\text{for $k=1,\ldots,N.$}
  \end{equation}
Now suppose $\uu\in C^2(\bOmg,\R^N)$ with $\uu=0$ on $\Gam$. Then,
after taking derivatives in $\Omg_s$ for some $s>0$ and then taking
the trace on $\Gam$, using \qref{E.n.id4.0} we get
\[
\div \( (\tanproj)\uu \)
= \pa_j \big( u_j - n_j n_k u_k \big)=\pa_j u_j -
n_j n_k \pa_j u_k =\pa_j u_j - \pa_k u_k =0.
\]
For general $\uu\in H^2(\Omg,\R^N)$ with $\uu|_\Gam=0$, the expression
$\div \( (\tanproj)\uu \)$
is in $H^1(\Omg_s)$ for small $s>0$ and hence is in $H^{1/2}(\Gamma)$
by a trace theorem. After approximating $\uu$ using Lemma~\ref{L.C2density}
we obtain the result in (i).

(ii) At first we suppose $\uu\in C^2(\bOmg,\R^N)$ with $\uu=0$ on $\Gam$.
We claim in fact that for any $f \in C^2(\bOmg)$ with
$f\big|_{\Gam}=0$,
  \begin{equation} \label{E.n.idf}
     \nn \cdot (\lap  - \grad \div )( \nn f )=0 \quad\text{on $\Gam$.}
  \end{equation}
This yields (ii) by taking $f = \nn \cdot \uu.$
We prove \qref{E.n.idf} in two steps.

1. The formula in (i) holds in $C(\Gamma)$ if $\uu$ is $C^1$.
Since $\tanproj=(\tanproj)^2$, we can use $\uu=(\tanproj)\grad f$ in
(i) to find that
\begin{equation}\label{E.L.f1}
   \div \( (\tanproj)\grad f \)= 0 \quad\text{ on $\Gam.$}
\end{equation}

2. Using \qref{E.n.id1} it is easy to verify
the following identities in $\Omg_s$:
\begin{eqnarray}
\nn\cdot\lap(\nn f) &=& \lap f + f\,\nn\cdot\lap\nn ,
\label{E.f.id1} \\
\nder\div(\nn f) &=&
(\nn\nn^t):\nabla^2f + (\div\nn)\nder f + f\nder\div\nn,
\label{E.f.id2} \\
\div(\nn\nn^t\grad f) &=&
(\nn\nn^t):\nabla^2f + (\div\nn)\nder f.
\end{eqnarray}
Here $(\nn\nn^t):\nabla^2f:=n_in_j\pa_i\pa_jf$. It directly follows that
  \begin{equation} \label{E.lemma1.3}
    \nn \cdot (\lap - \grad \div ) (\nn f) =
    \grad \cdot (\tanproj)\grad f + f \nn \cdot (\lap - \grad \div )\nn .
  \end{equation}
Using this with \qref{E.L.f1} proves \qref{E.n.idf},
and establishes (ii) when $\uu\in C^2(\bOmg)$ with $\uu=0$ on $\Gam$.

To establish (ii) for arbitrary $\uu\in H^2(\Omg,\R^N)$, we
restrict to $\Omg_s$ for small $s$ and let
$\ba=(\lap-\grad\div)(\nn\nn^t\uu)$. Then $\ba\in
L^2(\Omg_s,\R^N)$ and $\div\ba=0$ in the sense of distributions,
so $\ba\in\Hdivs$ and a well-known trace theorem (see \cite{GR},
theorem 2.5) yields that the map $H^2(\Omg_s,\R^N)\to\Hdivs\to
H^{-1/2}(\Gamma)$ given by $\uu\mapsto\ba\mapsto \nn\cdot\ba$ is
continuous. To conclude the proof, simply apply the approximation
lemma above to infer $\nn\cdot\ba|_\Gam=0$. $\square$

\subsection{Identities for the Stokes pressure} \label{S.upeupa}

Given $\uu \in \uspace$,
recall that $\P(\grad\div\uu)=0$, so that the Stokes pressure 
defined in \qref{StokesP} satisfies
\begin{equation}\label{pdef}
\grad \ps = \lap\uu - \grad \div \uu -\P\lap\uu = 
(I-\P)(\lap-\grad\div)\uu .
\end{equation}
Also recall that 
whenever $\ba\in L^2(\Omg,\R^N)$ and $\div\ba\in L^2(\Omg)$,
$\ndot\ba\in H^{-1/2}(\Gam)$ by the trace theorem for $\Hdiv$.
If $\div\ba=0$ and $\nn\cdot\ba|_\Gam=0$, then we have
$\<\ba,\grad\phi\>=0$ for all $\phi\in H^1(\Omg)$ 
and this means $(I-\P)\ba=0$. 
Thus, the Stokes pressure is not affected by any
part of the velocity field that contributes nothing to
$\nn\cdot\ba|_\Gam$ where $\ba=(\lap-\grad\div)\uu$.
Indeed, this means that the Stokes pressure is not affected by the
part of the velocity field in the interior of $\Omg$ away from the
boundary, nor is it affected by the normal component of velocity
near the boundary, since $\nn\cdot(\lap-\grad\div)(\nproj\uu)|_\Gam=0$
by Lemma~\ref{L.uH2}.

This motivates us to focus on the part of velocity near and parallel
to the boundary. We make the following decomposition.
Let $\rho:[0,\infty)\to[0,1]$ be a smooth decreasing function with
$\rho(t)=1$ for $t<\frac12$ and $\rho(t)=0$ for $t\ge1$.  For small
$s>0$, the cutoff function given by $\cut(x)=\rho(\Phi(x)/s)$
is $C^3$, with $\cut=1$ when $\Phi(x)<\frac12s$ and $\cut=0$
when $\Phi(x)\ge s$.  Then we can write
\begin{equation}\label{upaupe}
\uu = \upe + \upa
\end{equation}
where
\begin{equation} \label{updefs}
     \upe = (1-\cut) \uu + \cut \nproj \uu , \qquad
     \upa = \cut (\tanproj)\uu  .
\end{equation}

Since $\upe = (\nproj) \uu$ in $\Omg_{s/2}$, with
$\ba_{\perp}=(\lap-\grad\div)\upe$ we have
\begin{equation}
\mbox{$\ba_{\perp}\in L^2(\Omg,\R^N)$,\qquad
$\div\ba_{\perp}=0$\quad and\quad $\nn\cdot\ba_{\perp}|_\Gam=0$}
\end{equation}
by Lemma~\ref{L.uH2}(ii).  Hence $\<\ba_\perp,\grad\phi\>=0$ for
all $\phi\in H^1(\Omg)$, that is,
\begin{equation}\label{upezero}
(I-\P) (\lap - \grad \div )\upe=0.
\end{equation}
Combining this with \qref{pdef} and \qref{upaupe}
proves part (i) of the following.
\begin{lemma} \label{L.pair}
 Let $\Omg\subset\R^N$ be a bounded domain with $C^3$ boundary, and
 let $\uu \in \uspace$.
Let $\ps$ and $\upa$ be defined as in
 \qref{pdef} and \qref{updefs} respectively. Then
 \begin{description}
   \item (i) The Stokes pressure is determined by $\upa$ according to
   the formula
 \begin{equation}\label{psupa}
 \grad \ps = (I-\P) (\lap - \grad \div) \upa .
 \end{equation}
   \item (ii) For any $q \in H^1(\Omg)$ that satisfies
   $\lap q =0$ in the sense of distributions,
             \begin{equation} \label{stokespair-1}
               \< \lap \upa -\grad \ps, \grad q \> = 0.
             \end{equation}
   \item (iii) In particular we can let $q=\ps$ in (ii), so
   $ \< \lap \upa-\grad\ps, \grad \ps \>  = 0$ and
\begin{equation} \label{stokespair-2}
\| \lap \upa \|^2 = \| \lap \upa - \grad \ps \|^2 + \| \grad \ps \|^2.
              \end{equation}
 \end{description}
\end{lemma}
{\bf Proof:} We already proved (i). For (ii), note
by Lemma~\ref{L.uH2}(i) we have
\begin{equation}
     \div \upa |_{\Gam} =0,
\end{equation}
so $\div\upa\in H^1_0(\Omg)$, thus
$\<\grad\div\upa,\grad q\>=-\<\div\upa,\lap q\>=0$. Now (i) entails
\begin{equation}
\<\grad\ps,\grad q\> =
\<\lap\upa,\grad q\>.
\end{equation}
This proves (ii), and
then (iii) follows by the $L^2$ orthogonality. $\square$
\subsection{Proof of Theorem~\ref{T.main}}

Let $\eps>0$ and $\beta=\frac23+\eps$. We fix $\beta_1<1$ such that 
$1+\eps_0:=\beta(1+\frac12\beta_1^2)>1$, and fix $s>0$ small so
Theorem~\ref{T.N2D} (ii) applies in $\Omg_s$ with this $\beta_1$. 
Let $\uu\in \uspace$ and define the
Stokes pressure $\grad\ps$ by \qref{StokesP} and the decomposition
$\uu=\upe+\upa$ as in the previous subsection. 
Then by part (iii) of Lemma~\ref{L.pair} we have
\begin{equation}\label{lapu}
\|\lap\uu\|^2 = \|\lap\upe\|^2+ 2\<\lap\upe,\lap\upa\> +
\|\lap\upa-\grad\ps\|^2+\|\grad\ps\|^2.
\end{equation}
We will establish the Theorem with the help of two further estimates.

\noindent{\bf Claim 1:} For any $\eps_1>0$, there exists a constant
$C_1>0$ independent of $\uu$ such that
\begin{equation}\label{E.claim2}
\<\lap\upe,\lap\upa\> \ge -\eps_1\|\lap\uu\|^2 - C_1\|\grad\uu\|^2.
\end{equation}

\noindent{\bf Claim 2:} For any $\eps_1>0$ there exists a constant 
$C_2$ independent of $\uu$ such that
\begin{equation}
\|\lap\upa-\grad\ps\|^2 \ge 
\frac{\beta_1^2}{2} \|\grad\ps\|^2 -\eps_1\|\lap\uu\|^2 -
C_2\|\grad\uu\|^2.
\end{equation}

\noind {\bf Proof of claim 1:} From the definitions in \qref{updefs},
we have
\begin{equation}\label{upests}
\lap\upe = \cut \nproj \lap \uu + (1-\cut)\lap \uu + R_1,
\qquad
\lap\upa =
    \cut (\tanproj) \lap \uu  + R_2,
\end{equation}
where $\|R_1\|+\|R_2\|\le C\|\grad\uu\|$ with $C$ independent of
$\uu$.
Since $\tanproj=(\tanproj)^2$,
\[
  \( \cut \nproj \lap \uu +(1-\cut)\lap \uu  \)\cdot
  \( \cut (\tanproj) \lap \uu \) =
  0+\cut (1-\cut) | (\tanproj) \lap \uu |^2 \geq 0.
\]
This means the leading term of $  \<\lap \upe ,\lap \upa\> $ is
non-negative. Using the inequality $|\<a,b\>|\le
(\eps_1/C)\|a\|^2+(4C/\eps_1)\|b\|^2$ and the bounds on $R_1$ and
$R_2$ to estimate the remaining terms, it is easy to obtain
\qref{E.claim2}.

\noindent{\bf Proof of claim 2:}
Recall that $\upa$ is supported in $\Omgs$, and note
\begin{equation}
\lap\upa = \cut(\tanproj)\lap\uu + R_3
\end{equation}
where $\|R_3\|\le C\|\grad\uu\|$.
Since $\nn\cdot(\tanproj)\lap\uu=0$ we find 
\begin{equation}
\|\nn\cdot\lap\upa\|_{\Omgs} \le C_2 \|\grad\uu\|
\label{upaest00}
\end{equation}
with $C_2>0$ independent of $\uu$. 
We use $|a+b|^2\ge (1-\eps_2)|b|^2-|a|^2/\eps_2$ to get
\begin{align}
 \|\lap\upa&-\grad\ps\|^2_\Omg \ge  %
 \int_{\Omg_s^c} |\grad\ps|^2 + 
 \int_\Omgs |\nn\cdot(\lap\upa-\grad\ps)|^2
 \nonumber\\ \ge &\ %
 \int_{\Omg_s^c} |\grad\ps|^2 + 
 (1-\eps_2) \int_\Omgs |\nn\cdot\grad\ps|^2 
 - \frac{1}{\eps_2} \int_\Omgs |\nn\cdot\lap\upa|^2 .
% \nonumber\\ \ge &\ %
% \pcon_1\|\grad\ps\|^2_\Omg - C_2 \|\grad\uu\|^2_\Omg.
\label{upest0}
\end{align}
Next we use part (ii) of Theorem~\ref{T.N2D} with $p_0=0$ and
with $\beta_1\int_\Omgs|\nder p|^2$ added to both sides,
together with Lemma~\ref{L.pL2} and Poincar\'e's inequality,
to deduce that
\begin{equation}
\frac{\beta_1}{2} \int_\Omgs|\grad\ps|^2 \le
\int_\Omgs|\nder\ps|^2 + \eps_1\int_\Omgs |\lap\uu|^2 
+ C\int_\Omgs |\grad\uu|^2.
\label{upest1}
\end{equation}
Taking $1-\eps_2=\beta_1$ and combining \qref{upaest00},
\qref{upest0} and \qref{upest1} establishes Claim 2.

Now we conclude the proof of the theorem. Combining the two claims
with \qref{lapu}, we get
\begin{equation}
(1+3\eps_1) \|\lap\uu\|^2 \ge \left(1+\frac{\beta_1^2}2\right)
\|\grad\ps\|^2 - (C_2+2C_1)\|\grad\uu\|^2.
\end{equation}
Multiplying by $\beta$ and 
taking $\eps_1>0$ so that $3\eps_1<\eps_0$ concludes the proof.
$\square$

%%%%%%%%%%%%%%%%%%%%%%%%%%%%%%%%%%%%%%%%%%%%%%%%%%%%%%%%
\subsection{The space of Stokes pressures}
\label{S.SP}
%%%%%%%%%%%%%%%%%%%%%%%%%%%%%%%%%%%%%%%%%%%%%%%%%%%%%%%%%%%%%%%%%
According to \qref{mean0}--\qref{ps-wk}, 
the space of Stokes pressures, 
obtainable via \qref{StokesP}
from velocity fields $\uu\in \uspace$,
can be characterized
as the space 
\begin{equation}
\SP := 
\{ p\in H^1(\Omg)/\R\mid
\text{$\lap p=0$ in $\Omg$ and $\nder p|_\Gam\in \SPb$}\},
\end{equation}
where $\SPb$ is the subspace of $H^{-1/2}(\Gam)$ given by
\begin{equation} \label{xdef}
\SPb := \{ f= \nn\cdot(\lap-\grad\div)\uu|_\Gam \mid
\text{ $\uu\in \uspace $}\}.
\end{equation}
The Stokes pressure $p$ with zero average 
is determined uniquely by $f=\nder p|_\Gam\in \SPb$, 
with $\|p\|_{H^1(\Omg)}\le C\|f\|_{H^{-1/2}(\Gam)}$ 
by the Lax-Milgram lemma.

The space $\SPb$ may be characterized as follows.

\begin{theorem} \label{T.trace}
 Assume $\Omg \subset \R^N $ is a bounded, connected domain and its 
  boundary $\Gam$ is of class $C^3$.
Denote the connected components of $\Gam$ by $\Gam_i$, $i=1,\ldots,m$.
 Then 
\[
 \SPb = \{ f\in H^{-1/2}(\Gam) \mid  
\int_{\Gam_i} f =0 \text{\ \ for $i=1,\ldots,m$}\},
\]
and moreover, the map $\uu\mapsto \nn\cdot(\lap-\grad\div)\uu|_\Gam$ from 
$\uspace$ to $\SPb$ admits a bounded right inverse.
\end{theorem}

\noindent{\bf Proof.} First we check the necessity of the integral 
conditions. Let $u\in \uspace$ and let
$f= \nn\cdot(\lap-\grad\div)\uu|_\Gam$. 
For each connected component $\Gam_i$ of $\Gam$,
there is an $s_i>0$ small enough and a smooth
cut-off function $\rho_i$ defined in $\Omg$ which satisfies
$\rho_i(x) =1$ when $\dist(x,\Gam_i) < s_i$ and $\rho_i(x) =0$
when dist$(x,\Gam_j) < s_i$ for all $j \neq i$. 
Let $\ba= (\lap-\grad\div)(\rho_i\uu)$. Then 
$\ba\in L^2(\Omg,\R^N)$ and $\div\ba=0$, so
\begin{equation} \label{Gami-mean0}
  \int_{\Gam_i} f = \int_\Gam \nn\cdot\ba = \int_\Omg \div\ba =0.
\end{equation}

Next, let $f\in H^{-1/2}(\Gam)$ with $\int_{\Gam_i}f=0$ for all $i$.
Treating each boundary component separately, we can then solve the problem
\begin{equation}\label{s.lap}
\Delta_\Gamma \psi = -f    \quad\text{on $\Gam$}, 
\qquad \int_{\Gam_i}\psi=0 \quad\text{for $i=1,\ldots,m$},
\end{equation}
where $\Delta_\Gam$ is the (positive) Laplace-Beltrami operator on $\Gam$.
Denote the mapping $ f \mapsto \psi$ by $T$. 
Then $T \colon H^{-1}(\Gam) \to H^1(\Gam)$ is bounded (\cite[theorem
1.71, theorem 4.7]{Au}, \cite[p. 306, Proposition 1.6]{Ta}).
Also
$T\colon L^{2}(\Gam) \to H^2(\Gam)$ is bounded, by elliptic regularity
theory \cite[p. 306, Proposition 1.6]{Ta}.
So, interpolation implies (see \cite[vol I, p.~37, Remark 7.6]{LM})
\begin{equation}
  \|\psi \|_{H^{3/2}(\Gam)} \le C \|f\|_{H^{-1/2}(\Gam)}.
\end{equation}
Now by an inverse trace theorem \cite[Theorem 6.109]{RR}, 
there exists a map $\psi\mapsto q\in H^3(\Omg)$ with 
\begin{equation}\label{qvals}
\text{$q=0$\ \ and\ \ $\nder q=\psi$ \ \ on $\Gam$,}
\qquad \|q\|_{H^3(\Omg)}\le C\|\psi \|_{H^{3/2}(\Gam)}.
\end{equation}
We may assume $q$ is supported in a small neighborhood of $\Gam$.
Define
\begin{equation}\label{udef}
\uu = (\tanproj) \grad q.
\end{equation}
Then $f\mapsto\uu$ is bounded from $\SPb$ to 
$\uspace$. We claim
\begin{equation}\label{uclaim}
\nn\cdot(\lap-\grad\div)\uu = f \quad\text{on $\Gam$}.
\end{equation}

The proof of this claim amounts to showing, by calculations similar to those 
in the proof of Lemma~\ref{L.uH2}, that the normal derivative $\nder$ and
normal projection $\nproj$ commute on the boundary with the tangential gradient
and divergence operators $(\tanproj)\grad$ and $\div(\tanproj)$ for the
functions involved.

First, since $\nn\cdot\uu=0$, by expanding $\lap(\nn\cdot\uu)$ we get
\begin{equation}\label{tan1}
\nn\cdot\lap\uu = -(\lap\nn)\cdot\uu -2\grad\nn:\grad\uu =0 
\quad\text{on $\Gam$},
\end{equation}
since for each $i$, $\grad n_i$ is tangential and 
$\grad u_i$ is normal to $\Gam$ --- indeed, using $\pa_jn_i=\pa_in_j$ 
and \qref{E.n.id1} and \qref{E.n.id4.0}, we have that 
\begin{equation}\label{tan2}
\grad\nn:\grad\uu = (\pa_jn_i)(\pa_j u_i) = (\pa_i n_j)(n_jn_k\pa_ku_i)=0
\quad\text{on $\Gam$}.
\end{equation}
Next we calculate in $\Omg$ that
\begin{equation}\label{tan3}
\nder\div\uu = \div(\nder\uu)-\grad\nn:\grad\uu.
\end{equation}
Note that $\nder(\nproj)=0$ by \qref{E.n.id1}, so 
$\nder$ commutes with $\tanproj$ in $\Omg$.
Then since $\uu=(\tanproj)\uu$ from \qref{udef} we get
\begin{equation}\label{tan4}
\nder\uu = (\tanproj)(\nder)\uu =
(\tanproj)(\nder)\grad q.
\end{equation}
Now
\begin{equation}\label{tan5}
(\nder)\grad q = \grad(\nder q) - \ba
\end{equation}
where 
\begin{equation}\label{tan6}
a_i = (\pa_in_j)(\pa_j q) = (\pa_j n_i)(\pa_j q)
\end{equation}
This quantity lies in $H^2(\Omg)$ and vanishes on $\Gam$ 
since $\grad q=(\nproj)\grad q$ 
on $\Gam$. (This can be proved by approximation using Lemma~\ref{L.C2density}.)
Using part (i) of Lemma~\ref{L.uH2}, we have that 
$\div(\tanproj)\ba=0$ on $\Gam$.  Combining \qref{tan1}--\qref{tan5} 
we conclude that 
\begin{equation}\label{tan7}
\nn\cdot(\lap -\grad\div)\uu = -\div(\tanproj)\grad(\nder q)
\quad\text{on $\Gam$}.
\end{equation}
But it is well known
%\footnote{[[Is it? It should be. Citation??
%Do we need to include a proof?]]}
that at any point $x$ where $\Phi(x)=r\in(0,s)$, for any smooth
function $\phi$ on $\Omg_s$,
\begin{equation}
\div(\tanproj)\grad \phi = \lap \phi - (\div\nn)(\nder \phi) -
(\nder)^2 \phi
= -\lap_{\Gam_r}(\phi|_{\Gam_r}).
\end{equation}
where $\lap_{\Gam_r}$ is the Laplace-Beltrami operator on $\Gam_r$.
So taking $r\to0$ we see that the right hand side of \qref{tan7} is exactly 
$-\Delta_\Gamma(\nder q|_\Gamma)$.
So by \qref{s.lap} and \qref{qvals} we have established the claim in
\qref{uclaim}.  This finishes the proof.  $\square$

\medskip \noindent{\bf Remark 1.}
Given a velocity field $\uu\in H^2\cap H^1_0(\Omg,\R^3)$,
the associated Stokes pressure is determined by the normal
component at the boundary of the curl of the vorticity
$\omega=\curl\uu$, which is a vector field in $H^1(\Omg,\R^3)$.
A question related to Theorem~\ref{T.trace} is whether the space
$\SPb$ of such boundary values $\nn\cdot\curl\omega$ is constrained in
any way, as compared to the space of boundary values
$\nn\cdot\curl\vv$ where $\vv\in H^1(\Omg,\R^3)$ is arbitrary.

The answer is no. In \cite[Appendix I, Proposition 1.3]{Te}, Temam
proves 
\begin{equation}\label{curlH1}
  \curl H^1(\Omg,\R^3) = 
\{ \gg \in L^2(\Omg,\R^3)\mid \div \gg = 0, \int_{\Gam_i}\ndot \gg =0
\text{\ $\forall i$} \}.
\end{equation}
Clearly $\SPb\subset \nn\cdot\curl H^1(\Omg,\R^3)$ by \qref{xdef}.
For the other direction, let $\vv\in H^1(\Omg,\R^3)$ be
arbitrary, and let $f = \ndot\curl\vv|_\Gam$. By \qref{curlH1} or
otherwise, $f\in H^{-1/2}(\Gam)$ and $\int_{\Gam_i}f=0$ for all
$i$, hence $f\in \SPb$. This shows that for $N=3$,
\begin{equation}\label{xtemam}
\SPb=\ndot\curl H^1(\Omg,\R^3).
\end{equation}

A related point is that for $N=3$, the space of Stokes
pressure gradients $\grad\SP$ can be characterized as 
{\em the space of simultaneous gradients and curls}. 
\begin{theorem}\label{T.spgrad}
Assume $\Omg \subset \R^3 $ is a bounded, connected domain and its 
boundary $\Gam$ is of class $C^3$. Then
\begin{equation}
\grad\SP = \grad H^1(\Omg) \cap \curl H^1(\Omg,\R^3).
\end{equation}
\end{theorem}
\noind{\bf Proof.}
Indeed, $\grad\SP\subset\curl H^1$ by \qref{curlH1} and
Theorem~\ref{T.trace}. On the other hand, if $\gg=\curl\vv=\grad
p$ then $\lap p=\div\gg=0$ and $\nder p|_\Gam\in \SPb$ by
\qref{curlH1} and Theorem~\ref{T.trace}, so $\grad p\in\grad \SP$.
$\square$

\medskip \noindent{\bf Remark 2.}
In the book \cite{Te} (see Theorem 1.5)
Temam establishes the orthogonal decomposition 
$L^2(\Omg,\R^N)=H \oplus H_1 \oplus H_2$, 
which means that for any $g \in L^2(\Omg,\R^N)$,
\begin{equation} \label{Temam.decomp}
  \gg = \PP \gg +  \grad q + \grad \lap^{-1} \div \gg,
\end{equation}
where $q$ satisfies $\lap q = 0$ and $\ndot \grad q|_\Gam = \ndot
(\gg - \grad  \lap^{-1} \div \gg ) $.
By contrast, we have shown
\begin{equation} \label{new.decomp}
  \gg = \PP \gg +  \grad p + \grad \div \lap^{-1} \gg
\end{equation}
where $p$ satisfies $\lap p = 0$ and $\ndot \grad p|_\Gam = \ndot
(\gg -  \grad  \div \lap^{-1} \gg ) $, i.e., $p$ is the Stokes
pressure associated with $\lap^{-1} \gg$.
Thus the map $\gg\mapsto\grad p - \grad q$ is the commutator 
$\grad\lap^{-1}\div - \grad\div\lap^{-1}$. The 
decomposition \qref{Temam.decomp} is orthogonal, 
and $q$  satisfies $\< \ndot \grad q, 1 \>_{\Gam} = 0$. 
In our decomposition \qref{new.decomp}, the gradient terms are not 
orthogonal, but the Stokes pressure term enjoys the bounds 
stated in Corollary~\ref{C.1}, and if $\Gam$ is not connected, it has
the extra property that $\< \ndot \grad p, 1 \>_{\Gam_i} = 0$ for
every $i$.

\section{Unconditional stability of time discretization with pressure explicit}
\label{S.stab}

In this section we exploit Theorem~\ref{T.main} to establish the
unconditional stability of a simple time discretization scheme for the
initial-boundary-value problem for \qref{NSE5}, our unconstrained
formulation of the Navier-Stokes equations. We focus here on the case of
two and three dimensions. In subsequent sections we
shall proceed to prove an existence and uniqueness theorem based on
this stability result.

%%%%%%%%%%%%%%%%%%%%%%%%%%%%%%%%%%%%%%%%%%%%%%
Let $\Omg$ be a bounded domain in $\R^N$ with boundary $\Gam$ of class $C^3$.
We consider the initial-boundary-value problem
 \begin{align}
 \pa_t \uu + \udotgrad \uu + \nabla \pe + \nu \nabla \ps
 = \nu \Delta \uu + \ff
  &\qquad (t>0,\ x\in\Omg), \label{newNSE1} \\
  \uu = 0
  &\qquad (t\ge0,\ x\in\Gam), \label{newNSE2}\\
  \uu=\uuin
  &\qquad (t=0,\ x\in\Omg). \label{newNSE3}
 \end{align}
We assume $\uuin \in H_0^1(\Omg, \R^N)$ and
$\ff \in L^2(0,T;L^2(\Omg, \R^N))$ for some given $T>0$.
As before, the Euler and Stokes pressures $\pe$ and $\ps$ are defined by
the relations
 \begin{align}
  \PP( \udotgrad \uu - \ff) = \udotgrad \uu - \ff + \nabla \pe, &
  \label{E.EulerP} \\[4pt]
  \PP(-\Delta \uu ) = -\Delta \uu +\nabla(\div\uu) + \nabla \ps. &
  \label{E.StokesP}
 \end{align}
%%%%%%%%%%%%%%%%%%%%%%%%%%%%%%%%%%%%%%%%%%%%%%

Theorem~\ref{T.main} tells us that the Stokes pressure can be
strictly controlled by the viscosity term. This allows us to treat
the pressure term explicitly, so that the update of pressure is
decoupled from that of velocity. This can make corresponding
fully discrete numerical schemes very efficient (see \cite{JL}).
Here, through Theorem~\ref{T.main}, we will prove that the following
spatially continuous time discretization scheme
has surprisingly good stability properties:
\begin{align}
  \frac{\uu^{n+1} - \uu^n }{\Dt} -  \nu \lap \uu^{n+1}
= \ff^n -\uu^n \cdot \grad \uu^n - \grad \pe^n - \nu \grad \ps^n ,
\label{fdiff1}  \\
    \grad \pe^n =  (I-\P) (\ff^n- \uu^n \cdot \grad \uu^n ) ,
\label{fdiff2-1}\\
    \grad \ps^n = (I-\P)\lap \uu^n - \grad (\div\uu^n) ,
\label{fdiff2-2}\\
  \uu^n\big|_{\Gam} = 0.      \label{fdiff3}
\end{align}
We set
\begin{equation}\label{fndef}
\ff^n = \frac{1}{\Dt} \int_{n\Dt}^{(n+1)\Dt} \ff(t)\, dt,
\end{equation}
and take $\uu^0\in \uspace$
to approximate $\uuin$ in $H^1_0(\Omg,\R^N)$.
It is evident that for all $n=0,1,2,\ldots$, given $\uu^n\in H^2\cap H^1_0$
one can determine $\grad\pe^n\in L^2$ and $\grad\ps^n\in L^2$
from \qref{fdiff2-1} and \qref{fdiff2-2}
and advance to time step $n+1$ by solving \qref{fdiff1} as an elliptic
boundary-value problem with Dirchlet boundary values to obtain
$\uu^{n+1}$.  

This simple scheme is related to one studied by Timmermans et al.
\cite{Ti96}. In the time-differencing scheme described in \cite{Ti96}
for the linear Stokes equation, the pressure $p^n=\pe^n+\nu\ps^n$ is
updated in nearly equivalent fashion,
if one omits the velocity correction step that imposes zero divergence,
and uses first-order time differences in (15) and (18) of \cite{Ti96}.  
%\footnote{[[to omit?]] Note that $q^{n+1}=p^{n+1}-p^n+\nu\div\uu^{n+1}$ 
%satisfies
%\begin{equation}
%\< \grad q^{n+1},\grad \phi\> = \ip{ 
%\frac{\uu^{n+1}-\uu^n}{\Dt}+\ff^{n+1}-\ff^n+\uu^n\cdot\grad\uu^n-
%\uu^{n+1}\cdot\grad\uu^{n+1}, \grad\phi}
%\end{equation}
%for all $\phi\in H^1(\Omg)$. This compares with (18) and (23) of \cite{Ti96}.} 
Also see \cite{Pe,GuS,JL}.

Let us begin making estimates --- our main result is stated as
Theorem~\ref{T.stab} below.  
Dot \qref{fdiff1} with $-\lap u^{n+1}$ and use \qref{fdiff2-1}
and $\|I-\PP\|\le1$ to obtain
\begin{align}
  \frac{1}{2\Dt} \Big( &\|\grad \uu^{n+1}\|^2  - \| \grad \uu^n\|^2 +
  \| \grad \uu^{n+1} - \grad \uu^n\|^2 \Big) + \nu \|\lap \uu^{n+1}\|^2 
 \nonumber\\
  & \leq \|\lap \uu^{n+1} \| \Big(
  2\| \ff^n - \uu^n \cdot \grad \uu^n \|
 + \nu \|\grad \ps^n \| \Big)
\nonumber\\&
\leq \frac{\eps_1}{2}\|\lap\uu^{n+1}\|^2 + \frac{2}{\eps_1}
\|\ff^n-\uu^n\cdot\nabla\uu^n\|^2
%\\& \quad
+\frac{\nu}{2} \( \|\lap\uu^{n+1}\|^2+\|\grad\ps^n\|^2\)
\label{E.dot1}
\end{align}
for any $\eps_1>0$. (This is not optimal for $\grad\pe^n$ but is convenient.)
This gives
\begin{align}
  \frac{1}{\Dt} \Big( & \| \grad \uu^{n+1}\|^2 - \| \grad \uu^n\|^2
% Not needed! + \| \grad \uu^{n+1} - \grad \uu^n\|^2 
  \Big) +
  (\nu-\eps_1) \|\lap \uu^{n+1}\|^2 \nonumber
 \\ & \le
\frac{8}{\eps_1}\left( \|\ff^n\|^2+\|\uu^n\cdot\nabla\uu^n\|^2\right)
+\nu \|\grad \ps^n \|^2 .
  \label{E.temp.9}
\end{align}
Fix any $\pcon$ with $\frac23<\pcon<1$. 
By Theorem~\ref{T.main} one has
\begin{equation} \label{E.cite.them1}
 \nu \|\grad \ps^{n}\|^2 \le \nu\pcon \|\lap \uu^{n} \|^2
  +\nu C_{\pcon} \|\grad \uu^{n} \|^2.
\end{equation}
Using this in \qref{E.temp.9},  one obtains
%with $\nu(1-\pcon) = (\nu - \eps_1) -( \nu \pcon-\eps_1) $
\begin{align}
  \frac{1}{\Dt} \Big(  \|\grad \uu^{n+1}\|^2  -& \|\grad \uu^n\|^2 
% Not needed! +  \| \grad \uu^{n+1}   -  \grad \uu^n\|^2 
  \Big)
   + (\nu-\eps_1) \( \|\lap \uu^{n+1}\|^2 - \|\lap \uu^n\|^2 \)
   \nonumber\\
   &   + ( \nu - \eps_1 -\nu\pcon ) \|\lap \uu^n \|^2
    \nonumber \\
    &  \le \frac{8}{\eps_1}
 \( \|\ff^n\|^2 + \| \uu^n \cdot \grad \uu^n \|^2 \)
    + \nu C_{\pcon} \|\nabla \uu^n \|^2 .
    \label{E.fd-mid}
\end{align}

At this point there are no remaining difficulties with controlling the pressure.
It remains only to use the viscosity to control the nonlinear term.
We focus on the physically most interesting
cases $N=2$ and 3. We make use of Ladyzhenskaya's inequalities
\cite{La}
\begin{align} \label{E.lady1}
\int_{\R^N} g^4 &\le
2\left(\int_{\R^N} g^2\right) \left(\int_{\R^N} |\grad g|^2\right)
\qquad\qquad (N=2) ,\\
\int_{\R^N} g^4 &\le
4\left(\int_{\R^N} g^2\right)^{1/2} \left(\int_{\R^N} |\grad g|^2\right)^{3/2}
\qquad (N=3) ,\label{E.lady2}
\end{align}
valid for $g\in H^1(\R^N)$ with $N=2$ and $3$ respectively,
together with the fact that the standard
bounded extension operator $H^1(\Omg)\to H^1(\R^N)$ is also bounded
in $L^2$ norm,
%from $L^2(\Omg)$ to $L^2(\R^N)$,
to infer that for all $g\in
H^1(\Omg)$,
\begin{align}\label{E.lady3}
\|g\|_{L^4}^2 &\le C \|g\|_{L^2} \|g\|_{H^1} \hspace{3.4cm}(N=2),\\
\|g\|_{L^3}^2 &\le \|g\|_{L^2}^{2/3}\|g\|_{L^4}^{4/3}\le
C \|g\|_{L^2} \|g\|_{H^1} \qquad (N=3).\label{E.lady4}
\end{align}
Using that $H^1(\Omg)$ embeds into $L^4$ and $L^6$,
these inequalities lead to the estimates
\begin{align}
 \int_\Omg |\uu^n\cdot\grad\uu^n|^2 \le %\nonumber \\
 \begin{cases}
\|\uu^n\|_{L^4}^2 \|\grad \uu^n\|_{L^4}^2
\le C\|u\|_{L^2}\|\grad\uu^n\|_{L^2}^2\|\grad\uu^n\|_{H^1}
& (N=2), \\[4pt]
\|\uu^n\|_{L^6}^2 \|\grad \uu^n\|_{L^3}^2
\le C\|\grad\uu^n\|_{L^2}^3\|\grad\uu^n\|_{H^1}
& (N=3).
\end{cases}\label{E.lady5}
\end{align}
By the elliptic regularity estimate
$\|\grad \uu\|_{H^1} \le \|\uu\|_{H^2}\le C \|\lap \uu\| $,
we conclude
\begin{align} \label{E.nlin}
 \| \uu^n \cdot \grad \uu^n \|^2
%& \leq \| \uu^n \|^2_{L^6} \|\grad \uu^n \|^2_{L^3} \leq
%C \|\grad \uu^n \|^2_{L^2}\|\grad \uu^n \|_{L^2} \| \grad \uu^n \|_{H^1}
%\nonumber\\ &
& \leq 
\begin{cases}
\eps_2  \|\lap \uu^n \|^2
+ {4C}{\eps_2^{-1}}\|\uu^n\|^2 \| \grad \uu^n \|^4
& (N=2), \\[4pt]
\eps_2  \|\lap \uu^n \|^2
+ {4C}{\eps_2^{-1}} \| \grad \uu^n \|^6
& (N=2 \mbox{ or } 3).
\end{cases}
\end{align}
for any $\eps_2>0$. Plug this into \qref{E.fd-mid} and take $\eps_1$,
$\eps_2 >0$ satisfying $\nu - \eps_1 >0$ and 
$\eps:=\nu - \eps_1-\nu\pcon  - 8 \eps_2 / \eps_1 > 0 $. We get
\begin{align}
 \frac{1}{\Dt} \( & \|\grad \uu^{n+1}\|^2  - \|\grad \uu^n\|^2 \)
  + (\nu -\eps_1 ) \( \|\lap \uu^{n+1}\|^2 - \| \lap \uu^n \|^2 \) 
%  \nonumber \\
%  & + \frac{1}{\Dt} \|\grad \uu^{n+1} - \grad \uu^n \|^2
  + \eps \| \lap \uu^n \|^2
   % + \frac{1}{\eps_1} \| \grad \pe^n \|^2
    \nonumber \\
 &\leq   \frac{8}{\eps_1} \|\ff^n\|^2
   + \frac{32C}{\eps_1 \eps_2} \|\grad \uu^{n} \|^6
   + \nu C_{\pcon} \|\grad \uu^{n}\|^2 . \label{fdest1}
\end{align}
A simple discrete Gronwall-type argument
leads to our main stability result:
\begin{theorem}\label{T.stab}
Let $\Omg$ be a bounded domain in $\R^N$ ($N=2$ or $3$) with $C^3$
boundary, and assume $\ff \in L^2(0,T;L^2(\Omg, \R^N))$ for some
given $T>0$ and $\uu^0 \in H_0^1(\Omg, \R^N)\cap H^2(\Omg,\R^N)$.
 Consider the time-discrete scheme \qref{fdiff1}-\qref{fndef}.
 Then there exist positive constants
 $T^*$ and $C_3$,
 such that whenever $n\Dt\le T^*$, we have
 \begin{align}
  \sup_{0\leq k \leq n }\|\grad \uu^k \|^2
  + \sum_{k=0}^n  \|\lap \uu^k \|^2 \Dt \le C_3 ,
  %+ \sum_{k=0}^n \( \|\lap \uu^k \|^2 + \|\grad \ps^k\|^2 \) \Dt \le C_3 ,
  \label{stab1} \\
%% This estimate follows from the other two by parts & Cauchy-Schwarz
%   \sum_{k=0}^{n-1} \| \grad (\uu^{k+1}-\uu^k)\|^2 \leq C_3 ,
%   \label{stab2} \\
\sum_{k=0}^{n-1} \left( \left\|
\frac{\uu^{k+1}-\uu^k}{\Dt}\right\|^2 + \|\uu^k\cdot\grad\uu^k\|^2
\right) \Dt \le C_3. \label{stab3}
%+\|\grad\pe^k\|^2\right) \Dt \le C_3. \label{stab3}
 \end{align}
The constants $T^*$ and $C_3$ depend only upon $\Omg$, $\nu$ and
\[
M_0:=\|\grad\uu^0\|^2+\nu\Dt\|\lap\uu^0\|^2+\int_0^T\|\ff\|^2.
\]
\end{theorem}

\noindent{\bf Proof:} Put
\begin{equation}
z_n = \| \grad \uu^n \|^2 + (\nu-\eps_1) \Dt \|\lap \uu^n\|^2 ,
%\label{zdef} \\
\quad w_n = %\frac{1}{\Dt}\|\grad \uu^{n+1}-\grad\uu^n\|^2+
\eps \|\lap\uu^{n}\|^2 , %+ \frac{1}{\eps_1}\|\grad\pe^n\|^2 ,
%\label{wdef}\\
\quad b_n = \|\ff^n\|^2,
\end{equation}
%\begin{align}
%z_n &= \| \grad \uu^n \|^2 + (\nu-\eps_1) \Dt \|\lap \uu^n\|^2 ,
%\label{zdef} \\
%w_n &= %\frac{1}{\Dt}\|\grad \uu^{n+1}-\grad\uu^n\|^2+
%\eps \|\lap\uu^{n}\|^2 , %+ \frac{1}{\eps_1}\|\grad\pe^n\|^2 ,
%\label{wdef}\\
%b_n &= \|\ff^n\|^2,
%\end{align}
and note that from \qref{fndef} we have that as long as $n\Dt\le T$,
\begin{equation}\label{fest}
\sum_{k=0}^{n-1}\|\ff_k\|^2\Dt \le \int_0^T |\ff(t)|^2\,dt .
\end{equation}
Then by \qref{fdest1},
\begin{equation}
z_{n+1} +w_n\Dt \le z_n+ C\Dt(b_n+ z_n + z_n^3) , \label{zeq1}
\end{equation}
where we have replaced 
$\max \{ 8/\eps_1, 32C / (\eps_1 \eps_2),  \nu C_{\pcon}  \} $ by $C$.
Summing from 0 to $n-1$ and using \qref{fest}
%and the bound $\|\grad\ps^0\|\le C\|\lap\uu^0\|$
yields
\begin{equation}
z_n+\sum_{k=0}^{n-1}w_k\Dt \le CM_0 + C\Dt\sum_{k=0}^{n-1}(z_k+z_k^3)
=: y_n.
\label{ydef}
\end{equation}
The quantities $y_n$ so defined increase with $n$ and satisfy
\begin{equation}
y_{n+1}-y_n = C\Dt(z_n+z_n^3) \le C\Dt(y_n+y_n^3).
\label{yeq1}
\end{equation}
Now set $F(y)=\ln(\sqrt{1+y^2}/y)$ so that
$F'(y)=-(y+y^3)^{-1}$. Then on $(0,\infty)$,
$F$ is positive, decreasing and convex, and we have
\begin{equation}
F(y_{n+1})-F(y_n) = F'(\xi_n)(y_{n+1}-y_n) \ge
- \frac{y_{n+1}-y_n}{y_n+y_n^3} \ge -C\Dt,
\end{equation}
whence
\begin{equation}
F(y_n)\ge F(y_0)- Cn\Dt= F(CM_0)-Cn\Dt.
\end{equation}
Choosing any $T^*>0$ so that $C_*:= F(CM_0)-C T^*>0$, we infer that
as long as $n\Dt\le T^*$ we have $y_n\le F^{-1}(C_*)$, and this
together with \qref{ydef}
%and $\|\grad \ps^k  \|^2 \leq C \|\lap \uu^k \|^2$
yields the stability estimate \qref{stab1}.

Now, using \qref{E.nlin} and elliptic regularity, we get from
\qref{stab1} that
\begin{equation}
\sum_{k=0}^n \|\uu^k\cdot\grad\uu^k\|^2\Dt \le
C\sum_{k=0}^n \|\grad\uu^k\|_{L^2}^2\|\grad\uu^k\|_{H^1}^2\Dt \le
C\sum _{k=0}^n \|\lap\uu^k\|^2\Dt
\le C.
\end{equation}
Then the difference equation \qref{fdiff1} yields
\begin{equation}
\sum_{k=0}^{n-1} \left\| \frac{\uu^{k+1}-\uu^k}{\Dt}\right\|^2\Dt \le C.
\end{equation}
This yields \qref{stab3} and finishes the proof of the Theorem.
$\square$
%[[Note the $k=n$ term is not controlled only because of $\ff^n$!]]

%\endinput

\section{Existence and uniqueness of strong solutions} \label{S.exist}

The stability estimates in Theorem~\ref{T.stab} lead directly to the following
existence and uniqueness theorem for strong solutions of the unconstrained
formulation \qref{NSE5} of the Navier-Stokes equations. 
Regarding the constrained Navier-Stokes equations there are of course
many previous works; see \cite{Am} for a recent comprehensive treatment.
For unconstrained formulations of the Navier-Stokes equations with a
variety of boundary conditions including the one considered in the
present paper, Grubb and Solonnikov \cite{GS1,GS2} lay out
a general existence theory in anisotropic Sobolev spaces
using a theory of pseudodifferential initial-boundary-value problems
developed by Grubb.

\begin{theorem}\label{existence}
Let $\Omg$ be a bounded domain in $\R^3$ with boundary $\Gam$ of
class $C^3$, and let $\ff \in L^2(0,T;L^2(\Omg, \R^N))$, $\uuin
\in H_0^1(\Omg, \R^N)$. Then, there exists $T^*>0$ depending only upon
$\Omg$, $\nu$ and $M_1:=\|\grad\uuin\|^2+\int_0^T\|\ff\|^2$, 
so that a unique strong solution of \qref{newNSE1}-\qref{newNSE3}
exists on $[0,T^*]$, with
\begin{align*}
& \uu\in
L^2(0,T^*;H^2(\Omg,\R^N))\cap H^1(0,T^*;L^2(\Omg,\R^N)),\\
& \grad p= \grad \pe+\nu \grad \ps\in L^2(0,T^*;L^2(\Omg,\R^N)) ,
\end{align*}
where $\pe$ and $\ps$ are as in \qref{EulerP} and
\qref{StokesP}. Moreover, $\uu\in C([0,T^*],H^1(\Omg,\R^N))$,
and $\div\uu \in C^{\infty}((0,T^*],C^\infty(\Omg))$ is a
classical solution of the heat equation with no-flux boundary
conditions. The map $t\mapsto\|\div\uu\|^2$ is smooth for $t>0$
and we have the dissipation identity
\begin{equation} \label{E.dissip.iden}
  \frac{d}{dt}\frac12\|\div \uu \|^2 +  \nu \|\grad (\div \uu)\|^2 =0.
\end{equation}
\end{theorem}

\noindent
{\bf Proof of existence:}
% of Theorem \ref{existence}} \label{sectionproofexistence}
We shall give a simple proof of existence based on the finite difference
scheme considered in section \ref{S.stab}, using a classical
compactness argument \cite{Ta1,Te,LM}.
However, in contrast to similar arguments in other sources,
for example by Temam \cite{Te} for a time-discrete scheme with
implicit differencing of pressure terms, we do not make any use of
regularity theory for stationary Stokes systems.

First we smooth the initial data. Given $\uuin\in
H^1_0(\Omg,\R^N)$ and $\Dt>0$, determine $\uu^0$ in
$H^1_0\cap H^2(\Omg,\R^N)$ by solving
$(I-\Dt\,\lap)\uu^0=\uuin$. An energy estimate yields
\[
\|\grad\uu^0\|^2+\Dt\|\lap\uu^0\|^2 \le \|\grad\uuin\|
\|\grad\uu^0\| \le \|\grad\uuin\|^2.
\]
Then $\|\Dt\,\lap\uu^0\|^2=O(\Dt)$ as $\Dt\to0$, so
$\uu^0\to\uuin$ strongly in $L^2$ and weakly in $H^1$. The
stability constant $C_3$ in Theorem~\ref{T.stab} is then uniformly
bounded independent of $\Dt$.

We define the discretized solution $\uu^n$ by \qref{fdiff1}-\qref{fdiff3}
of section \ref{S.stab}, and note
\begin{equation}\label{udiff1}
  \frac{\uu^{n+1} - \uu^n }{\Dt} +
  \P( \uu^n \cdot \grad \uu^n -\ff^n -\nu\lap\uu^n)
     = \nu\lap(\uu^{n+1}-\uu^n)+ \nu \grad\div \uu^n .
\end{equation}

With $t_n=n\Dt$,
we put $\uuD(t_n)=\uud(t_n)=\uu_n$ for $n=0,1,2,\ldots$, and
define $\uuD(t)$ and $\uud(t)$ on each subinterval $[t_n,t_n+\Dt)$
through linear interpolation and as piecewise constant respectively:
\begin{align}  \label{E.uuD}
  \uuD(t_n+s)  &=
  \uu^n + s \left(\frac{\uu^{n+1}-\uu^n}{\Dt}\right),
%  (1-\frac{s}{\Dt}) \uu^n + \frac{s}{\Dt} \uu^{n+1},
  \qquad  s \in [0,\Dt), \\
  \uud(t_n+s) &= \uu^n, \hspace{4cm}  s \in [0,\Dt).
  \label{E.uud}
\end{align}
Then \qref{udiff1} means that whenever $t>0$ with $t\ne t_n$,
\begin{equation} \label{udiff2}
\D_t\uuD + \P( \uud\cdot\grad\uud -\ffd -\nu\lap\uud)
=\nu\lap(\uud(\cdot+\Dt)-\uud)+\nu\grad\div\uud,
\end{equation}
where $\ffd(t)=\ff^n$ for $t\in[t_n,t_n+\Dt)$.

We will use the simplified notation $X(Y)$ to
denote a function space of the form $X([0,T^*],Y(\Omg,\R^N))$,
and we let $Q=\Omg\times[0,T^*]$ where $T^*$ is given by
Theorem~\ref{T.stab}.
The estimates in Theorem~\ref{T.stab} say
that $\uuD$ is bounded in the Hilbert space
\begin{equation}\label{Vdef}
V_0:= L^2(H^2\cap H^1_0)\cap H^1(L^2),
\end{equation}
and also that $\uud$ is bounded in $L^2(H^2)$,
uniformly for $\Dt>0$.
Moreover, estimate \qref{stab1} says $\uuD$ is bounded in
$C(H^1)$. This is also a consequence of the embedding
$V_0\emb C(H^1)$, see \cite[p.\ 42]{Ta1} or \cite[p.\ 288]{Ev}.

Along some subsequence $\Dt_j\to0$, then, we have that $\uuD$
converges weakly in $V_0$ to some $\uu\in V_0$, and $\uud$ and
$\uud(\cdot+\Dt)$ converge
weakly in $L^2(H^2)$ to some $\UU_1$ and $\UU_2$ respectively.
Since clearly $V_0\emb H^1(Q)$, and since the embedding
$H^1(Q)\emb L^2(Q)$ is compact,
we have that $\uuD\to\uu$ strongly in $L^2(Q)$.
Note that by estimate \qref{stab3},
\begin{equation}\label{Dd1}
\|\uuD-\uud\|^2_{L^2(Q)} \le \|\uud(\cdot+\Dt)-\uud\|^2_{L^2(Q)} =
\sum_{k=0}^{n-1}\|\uu^{n+1}-\uu^n\|^2\Dt \le C\Dt^2.
\end{equation}
Therefore $\uud(\cdot+\Dt)$ and $\uud$ converge to $\uu$ strongly
in $L^2(Q)$ also, so $\UU_1=\UU_2=\uu$.

We want to show $\uu$ is a strong solution of \qref{newNSE1}
by passing to the limit in \qref{udiff2}.
From the definition of $\ff^n$ in \qref{fndef},
it is a standard result
which can be proved by using a density argument that
\[
\|\ff-\ffd\|_{L^2(Q)}^2 \to 0 \quad\mbox{as $\Dt\to0$.}
\]
We are now justified in passing to the limit weakly in $L^2(Q)$ in all
terms in \qref{udiff2} except the nonlinear term,
which (therefore) converges weakly to some $\vec{w}\in L^2(Q)$.
But since $\grad\uud$ converges
to $\grad\uu$ weakly and $\uud$ to $\uu$ strongly in $L^2(Q)$,
we can conclude $\uud\cdot\grad\uud$ converges to $\uu\cdot\grad\uu$
in the sense of distributions on $Q$. So $\vec{w}=\uu\cdot\grad\uu$, and
upon taking limits in \qref{udiff2} it follows that
\begin{equation}\label{E.again}
\pa_t \uu + \P \( \uu \cdot \grad \uu - \ff - \nu \lap
  \uu\) = \nu \grad \div \uu.
\end{equation}
That is, $\uu$ is indeed a strong solution of \qref{newNSE1}.
That $\uu(0)=\uuin$ is a consequence of the continuity of the map
$\uu\to\uu(0)$ from $V_0$ through $C(H^1)$ to $H^1(\Omg,\R^N)$.

It remains to study $\div\uu$.
Dot \qref{E.again} with $\grad\phi$, $\phi\in H^1(\Omg)$.
We get
\begin{equation}
  \int_\Omg \pa_t \uu \cdot\nabla \phi =
  \nu \int_\Omg \nabla(\div\uu) \cdot \nabla \phi.
  \label{heat3}
\end{equation}
This says that $w=\div\uu$ is a weak solution of the heat equation
with Neumann boundary conditions:
\begin{equation}
\pa_t w=\nu\lap w \quad\text{in $\Omg$},
\qquad \nder w=0 \quad\text{on $\Gam$}.
\end{equation}
Indeed, the operator $A:=\nu\lap$ defined on $L^2(\Omg)$ with domain
\begin{equation} \label{domA}
D(A) = \{ w\in H^2(\Omg)\mid \nn\cdot\grad w = 0 \mbox{ on $\Gamma$}\}
\end{equation}
is self-adjoint and non-positive, so generates an analytic semigroup.
For any $\phi\in D(A)$ we have that
$t\mapsto\<w(t),\phi\>=-\<u(t),\grad\phi\>$ is absolutely continuous,
and using \qref{heat3} we get
$(d/dt)\<w(t),\phi\> = \<w(t),A\phi\>$ for a.\,e.\ $t$.
By Ball's characterization of weak solutions of abstract
evolution equations \cite{B}, $w(t)=e^{At}w(0)$ for all $t\in[0,T^*]$.
It follows $w\in C([0,T^*], L^2(\Omg))$, and 
$w(t)\in D(A^m)$ for every $m>0$ \cite[theorem 6.13]{Pa}.
Since $A^mw(t)=e^{A(t-\tau)}A^m w(\tau)$ if $0<\tau<t$ we infer
that for $0<t\le T^*$, $w(t)$ is analytic in $t$ with values in $D(A^m)$.
Using interior estimates for elliptic equations, we find
$w\in C^\infty((0,T^*],C^\infty(\Omg))$ as desired.
The dissipation identity follows by dotting with $w$.

This finishes the proof of existence. $\square$

\pagebreak %\medskip
\noindent{\bf Proof of uniqueness:}
Suppose $\uu_1$ and $\uu_2$ are both solutions of
\qref{newNSE1}--\qref{newNSE3} belonging to $V_0$.
Put $\uu = \uu_1 -\uu_2$ and $\grad\ps=(I-\PP)(\lap-\grad\div)\uu$.  Then
$\uu(0)=0$ and
\begin{equation} \label{proofuniq1}
 \pa_t \uu + \P \(
% \uu_1 \cdot \grad \uu_1 -  \uu_2 \cdot \grad \uu_2
 \uu_1 \cdot \grad \uu +  \uu \cdot \grad \uu_2 )
= \nu \lap \uu -\nu\grad\ps .
\end{equation}
Dot with $-\lap \uu$ and use
Theorem~\ref{T.main} to get
\begin{equation}
\<\nu\lap\uu-\nu\grad\ps,-\lap\uu\> \le
-\frac{\nu}{2}\|\lap\uu\|^2+\frac{\nu}2\|\grad\ps\|^2
\le -\frac{\nu\pcon}{2}\|\lap\uu\|^2+ C\|\grad\uu\|^2.
\end{equation}
Next, use the Cauchy-Schwarz inequality for the nonlinear terms,
estimating them as
follows in a manner similar to \qref{E.lady3}-\qref{E.lady5}, using that
$\uu_1$ and $\uu_2$ are a priori bounded in $H^1$ norm:
\begin{align}
\| \uu_1 \cdot \grad \uu \|\, \|\lap\uu\|
&\leq C\|\grad\uu_1\| \|\grad\uu\|^{1/2} \|\lap\uu\|^{3/2}
\leq \eps \|\lap\uu\|^2 + C\|\grad\uu\|^2,
\\[3pt]
\| \uu \cdot \grad \uu_2 \|\, \|\lap\uu\|
&\leq C \|\grad\uu\| \|\grad\uu_2\|_{H^1} \|\lap \uu\|
\leq \eps \|\lap\uu\|^2 + C\|\lap\uu_2\|^2\|\grad\uu\|^2.
\end{align}
Lastly, since $\uu\in V_0$ we infer that
$\<\pa_t\uu,-\lap\uu\> \in L^1(0,T)$ and
$t\mapsto\|\grad\uu\|^2$ is absolutely continuous with
\begin{equation}\label{E.diss2}
\<\pa_t\uu,-\lap\uu\> =
\frac12 \frac{d}{dt} \|\grad\uu\|^2 .
\end{equation}
%See \qref{E.diss.1} for a remark on the proof of \qref{E.diss2}.
This can be shown by using the density of smooth functions in $V_0$; see
\cite[p.\ 287]{Ev} for a detailed proof of a similar result.

Through this quite standard-style approach, we get
\begin{equation}
  \frac{d}{dt} \|\grad \uu\|^2 + \alpha \|\lap \uu \|^2 \leq
  C (1+\|\lap \uu_2\|^2 ) \| \grad \uu \|^2
\end{equation}
for some positive constants $\alpha $ and $C$. Because $ \|\lap
\uu_2\|^2 \in L^1(0,T),$ by Gronwall's inequality we get $\|\grad
\uu \|\equiv0.$ This proves the uniqueness. $\square$

Since the interval of existence $[0,T_*]$ depends only upon
$M_1$, in standard fashion we may extend the unique strong solution
to a maximal interval of time, and infer that the approximations
considered above converge to this solution up to the maximal time.

\begin{corollary} Given the assumptions of Theorem~\ref{existence},
system \qref{newNSE1}-\qref{newNSE3} 
admits a unique strong solution $\uu$ on a maximal interval
$[0,\Tmax)$ with the property that if $\Tmax<T$ then 
\begin{equation}\label{blowup}
\|\uu(t)\|_{H^1}\to\infty \quad\mbox{as $t\to\Tmax$.}
\end{equation}
For every $\hat T\in[0,\Tmax)$,
the approximations $\uu_{\Dt}$ constructed in \qref{E.uuD} converge to
$\uu$ weakly in 
\[
L^2([0,\hat T],\uspace)\cap H^1([0,\hat T], L^2(\Omg,\R^N))
\]
and strongly in $L^2([0,\hat T]\times\Omg,\R^N)$.
\end{corollary}

\section{Unconditional stability and convergence for $C^1/C^0$
finite element methods without inf-sup conditions}
\label{S.fem}

The simplicity of the stability proof for the time-discrete scheme in
section 4 allows us to easily establish the unconditional stability 
and convergence (up to the maximal time of existence for the strong
solution) of
corresponding fully discrete finite-element methods that use $C^1$
elements for the velocity field and $C^0$ elements for pressure. 
To motivate the discretization, we write the unconstrained
Navier-Stokes formulation \qref{NSE6} in weak form as follows,
in terms of total pressure $p=\pe+\nu\ps$:
\begin{align}
\<\uu_t+\grad p -\nu\lap\uu +\udotgrad\uu-\ff&,\lap \vv\>=0
\quad\forall v\in\uspace,
\label{nse-var1} \\
\hspace{-0.3cm}
\<\grad p +\nu\grad\div\uu-\nu\lap\uu +\udotgrad\uu-\ff&,\grad \phi\>=0
\quad\forall \phi\in H^1(\Omg).
\label{nse-var2}
\end{align}

We suppose that for some sequence of positive values of $h$
approaching zero, $\xuh\subset\uspace$ is a finite-dimensional space
containing the approximate velocity field, and suppose $\xph\subset
H^1(\Omg)/\R$ is a finite-dimensional space containing approximate
pressures. We assume these spaces have the approximation property that
\begin{align}
\forall \vv\in\uspace\ &\forall h \ \exists\vv_h\in\xuh,
\quad \|\lap(\vv-\vv_h)\|\to 0 \quad\mbox{as $h\to0$},
\label{E.appx1}\\
\forall \phi\in H^1(\Omg)/\R\ &\forall h\  \exists \phi_h\in\xph,
\quad \|\grad(\phi-\phi_h)\|\to0 \quad\mbox{as $h\to0$}.
\label{E.appx2}
\end{align}
As we have emphasized in the introduction to this paper, we impose
{\em no} inf-sup condition between the spaces
$\xuh$ and $\xph$.
(We remark that in general, practical finite element methods usually use spaces
defined on domains that approximate the given $\Omg$.  For simplicity
here we suppose $\Omg$ can be kept fixed, such that finite-element
spaces $\xuh$ and $\xph$ can be found as described with $C^1$ elements for
velocity and $C^0$ elements for pressure.  Though generally
impractical, in principle this should be possible whenever $\Omg$ has
a piecewise polynomial $C^3$ boundary.) 

We discretize \qref{nse-var1}-\qref{nse-var2} in a straightforward
way, implicitly only in the viscosity term and explicitly in the
pressure and nonlinear terms. The resulting scheme was also derived 
in \cite{JL} and is equivalent to a
space discretization of the scheme in \qref{fdiff1}--\qref{fndef}.
Given the approximate velocity
$\uu^h_n$ at the $n$-th time step, we determine 
$p_h^n\in\xph$ and $\uu_h^{n+1}\in\xuh$ by requiring
\begin{align} 
\<\grad p_h^n +\nu\grad\div\uu_h^n-\nu\lap\uu_h^n
+\uu_h^n\!\cdot\!\grad\uu_h^n-\ff^n,\grad \phi_h\>=0
%\\
% & \< \grad p_h^n, \grad \phi_h\> = 
% \< \curl \uu_h^n , \nn \times \grad \phi_h\>_\Gam 
\quad &\forall \phi_h \in \xph, \label{fem.1}
  \\
  \< \frac{\grad\uu^{n+1}_h - \grad\uu^n_h }{\Dt}, \grad \vv_h \>  + 
\<\nu\lap\uu_h^{n+1},\lap\vv_h\>
%\nonumber \\ & \qquad 
= \<\grad p_h^n 
+\uu_h^n\!\cdot\!\grad\uu_h^n & -\ff^n,\lap\vv_h \>
%  \< \lap\uu^{n+1}_h, \lap \vv_h \> = \< \grad p_h^n,  \lap \vv_h \> = 
%  - \< \ff_h^n, \lap \vv_h \> 
\nonumber \\ 
  & \forall \vv_h \in \xuh . 
\label{fem.2} 
\end{align}

{\bf Stability.}
We are to show the scheme above is unconditionally stable. 
First, we take $\phi_h = p_h$ in \qref{fem.1}. Due to the fact that
\[
\<\PP(\lap-\grad\div)\uu_h^n,\grad p_h^n\>=0,
\]
we directly deduce from the Cauchy-Schwarz inequality that
\begin{equation}\label{E.phest}
  \|\grad p_h^n \| \le
  \|\nu\grad\ps(u_h^n)\|+\|\uu_h^n\!\cdot\!\grad\uu_h^n - \ff^n\|
\end{equation}
where
\begin{equation}  
  \grad\ps(u_h^n) = (I-\PP)(\lap-\grad\div)\uu_h^n
\end{equation}
is the Stokes pressure associated with $\uu_h^n$. (Note
$\grad\ps(u^n_h)$ need not lie in the space $\xph$). 
Now, taking $\vv_h = \uu_h^{n+1}$ in \qref{fem.2} and arguing just
as in \qref{E.dot1}, 
we obtain an exact analog of \qref{E.temp.9}, namely
\begin{align}
  \frac{1}{\Dt} \Big( & \| \grad \uu_h^{n+1}\|^2 - \| \grad \uu_h^n\|^2
% Not needed! + \| \grad \uu^{n+1} - \grad \uu^n\|^2 
  \Big) +
  (\nu-\eps_1) \|\lap \uu_h^{n+1}\|^2 \nonumber
 \\ & \le
\frac{8}{\eps_1}\left( \|\ff^n\|^2+\|\uu_h^n\cdot\nabla\uu_h^n\|^2\right)
+\nu \|\grad \ps(\uu_h^n) \|^2 .
  \label{fem.ineq}
\end{align}
Proceeding now exactly as in section 4 leads to the following
unconditional stability result. 

\begin{theorem}\label{T.stabh}
Let $\Omg$ be a bounded domain in $\R^N$ ($N=2$ or $3$) with $C^3$
boundary, and suppose spaces $\xuh\subset\uspace$, $\xph\subset
H^1(\Omg)/\R$ satisfy \qref{E.appx1}--\qref{E.appx2}.
Assume $\ff \in L^2(0,T;L^2(\Omg, \R^N))$ for some
given $T>0$ and $\uu_h^0 \in \xuh$.
 Consider the finite-element scheme \qref{fem.1}-\qref{fem.2} with
 \qref{fndef}.  Then there exist positive constants
 $T^*$ and $C_4$,
 such that whenever $n\Dt\le T^*$, we have
 \begin{align}
  \sup_{0\leq k \leq n }\|\grad \uu^k_h \|^2
  + \sum_{k=0}^n  \|\lap \uu^k_h \|^2 \Dt \le C_4 ,
  %+ \sum_{k=0}^n \( \|\lap \uu^k \|^2 + \|\grad \ps^k\|^2 \) \Dt \le C_3 ,
  \label{stab1h} \\
%% This estimate follows from the other two by parts & Cauchy-Schwarz
%   \sum_{k=0}^{n-1} \| \grad (\uu^{k+1}-\uu^k)\|^2 \leq C_3 ,
%   \label{stab2} \\
\sum_{k=0}^{n-1} \left( \left\|
\frac{\uu^{k+1}_h-\uu^k_h}{\Dt}\right\|^2 + \|\uu^k_h\cdot\grad\uu^k_h\|^2
\right) \Dt \le C_4. \label{stab3h}
%+\|\grad\pe^k\|^2\right) \Dt \le C_3. \label{stab3}
 \end{align}
The constants $T^*$ and $C_4$ depend only upon $\Omg$, $\nu$ and 
\[
M_{0h}:=\|\grad\uu_h^0\|^2+\nu\Dt\|\lap\uu_h^0\|^2+\int_0^T\|\ff\|^2.
\]
\end{theorem}

{\bf Convergence.} We prove the convergence of the finite-element
scheme described above by taking $h\to0$ to obtain the solution
of the time-discrete scheme studied in section, then $\Dt\to0$ as
before. Because of the uniqueness of the solution of the time-discrete
scheme and of the strong solution of the PDE, it suffices to prove
convergence for some subsequence of any given sequence of values of
$h$ tending toward $0$.  The bounds obtained in Theorem~\ref{T.stabh}
make this rather straightforward.

Fix $\Dt>0$. The bounds in Theorem~\ref{T.stabh} and in \qref{E.phest}
imply that for all positive integers $n<T_*/\Dt$, the $\uu^n_h$
are bounded in $\uspace$ and the $\grad p^n_h$ are bounded in
$L^2(\Omg,\R^N)$ uniformly in $h$. So from any sequence of $h$
approaching zero, we may extract a subsequence along which we have
weak limits
\begin{equation}\label{E.wklim}
 \uu^n_h \rightharpoonup \uu^n 
\mbox{ in $H^2(\Omg,\R^N)$},
\quad \grad p^n_h \rightharpoonup \grad p^n,
\ \uu^n_h\cdot\grad\uu^n_h \rightharpoonup \ww^n
\mbox{ in $L^2(\Omg,\R^N)$}
\end{equation}
for all $n$. Then $\uu^n_h\to\uu^n$ strongly in $H^1_0(\Omg,\R^N)$ and
so $\ww^n=\uu^n\cdot\grad\uu^n$ since the nonlinear term converges 
strongly in $L^1$.

Now, for any $\vv\in\uspace$ and $\phi\in H^1(\Omg)$, by assumption
there exist $\vv_h\in\xuh$, $\phi_h\in H^1(\Omg)$ such that 
$\vv_h\to\vv$ strongly in $H^2(\Omg,\R^N)$ and
$\grad\phi_h\to\grad\phi$ strongly in $L^2(\Omg,\R^N)$.
Applying these convergence properties in \qref{fem.1}--\qref{fem.2} 
yields that the weak limits in \qref{E.wklim} satisfy
\begin{align}
&\<\grad p^n +\nu\grad\div\uu^n-\nu\lap\uu^n
+\uu^n\!\cdot\!\grad\uu^n-\ff^n,\grad \phi\>=0 ,
\label{E.plim}
\\ &
\< \frac{\uu^{n+1} - \uu^n }{\Dt}
- \nu\lap\uu^{n+1}
+ \grad p^n +\uu^n\!\cdot\!\grad\uu^n  -\ff^n,\lap\vv \> = 0.
\label{E.ulim}
\end{align}
But this means exactly that $\uu^n$ satisfies \qref{fdiff1}
with $p^n=\pe^n+\nu\ps^n$, where
$\pe^n$ and $\ps^n$ are given by \qref{fdiff2-1}--\qref{fdiff2-2}.
So in the limit $h\to0$ we obtain the solution of the time-discrete
scheme studied in section~\qref{S.stab}. Then the limit $\Dt\to0$
yields the unique strong solution on a maximal time interval 
as established in section~\ref{S.exist}.

\section{Semigroup approach for the homogeneous linear case}
\label{S.semig}

There are many other approaches to existence theory for the
Navier-Stokes equations, of course --- Galerkin's method, mollification,
semigroup theory, etc. We will not discuss any of them here,
except to note that the linearization of the unconstrained
system \qref{NSE6} can be treated easily by analytic semigroup theory
using Theorem~\ref{T.main}.
Take $\nu=1$ without loss of generality, and
consider \qref{NSE6} without the nonlinear and forcing terms, i.e.,
consider the unconstrained Stokes equation
\begin{equation}\label{E.stokes0}
  \uu_t - \lap \uu + \grad\ps =0\qquad (t>0,\ x\in\Omg),
\end{equation}
with the no-slip boundary condition \qref{newNSE2} and initial
condition \qref{newNSE3}, where $\grad\ps$ is given by
\qref{StokesP} as before.
In the space $X=L^2(\Omg,\R^N)$ define operators $B_0$ and $B_1$ by
\begin{equation}\label{opdef}
B_0\uu = -\lap \uu, \qquad
B_1\uu = \grad\ps = (I-\PP)\lap \uu - \grad \div \uu ,
\end{equation}
with domain $D(B_0)=D(B_1)= \uspace$.
Then $B_0$ is a positive self-adjoint operator in $X$ with compact resolvent,
and by using
Theorem~\ref{T.main} together with the interpolation estimate
\[
\|\grad\uu\| \le \eps \|\lap\uu\| + C_\eps \|u\|
\]
valid for any $\eps>0$ for all $\uu\in D(B_0)$,
we deduce that
\begin{equation}\label{b1est}
\|B_1\uu\| \le a \|B_0\uu\| + K\|\uu\|
\end{equation}
for all $\uu\in D(B_0)$, where $a$ and $K$ are positive constants,
with $a<1$.
\begin{theorem} \label{T.B}
The unconstrained Stokes operator $B=B_0+B_1$ in the space 
$X=L^2(\Omg,\R^N)$ is sectorial and generates an analytic semigroup.
The resolvent of $B$ is compact and the spectrum of $B$
consists entirely of isolated eigenvalues of finite multiplicity, 
all of which are positive.  
Moreover, for any $\alpha\ge0$, given $\uuin\in D(B^\alpha)$ 
equation \qref{E.stokes0} has the solution
\[
\uu=e^{-Bt}\uuin\in C([0,T],D(B^\alpha))\cap C^\infty((0,T],D(B^m))
\]
for any $T>0$ and all $m>0$,
and this is the unique weak solution of $\pa_t\uu+B\uu=0$,
$\uu(0)=\uuin$ in the sense of Ball~\cite{B}.
\end{theorem}

\noindent{\bf Proof.}
That $B$ is sectorial is a consequence of \qref{b1est} and the
self-adjointness of $B_0$. Indeed, by a theorem on the perturbation
of sectorial operators \cite[p. 19, theorem 1.3.2]{He},
it suffices to show that for some
$\phi_0<\pi/2$,
\begin{equation}\label{Sest}
a \sup_{\lambda\in S_0} \|B_0(\lambda-B_0)^{-1}\| <1
\end{equation}
where $S_0\subset\C$ is the sector where
$\phi_0<|\arg\lambda|\le\pi$.
By expanding any element of (complexified) $X$ with respect to an
orthonormal basis of eigenfunctions of $B_0$,
for any $\lambda\notin\sigma(B_0)$ we get
% not an eigenvalue of $B_0$ we have
\[
\|B_0(\lambda-B_0)^{-1}\| = \sup_{\mu\in\sigma(B_0)}
\left|\frac{\mu}{\lambda-\mu}\right|
.
\]
Fix $\tilde a\in(a,1)$.
For any $\mu>0$, we have $|\mu|\le|\lambda-\mu|$ whenever $\Re\lambda\le0$,
and it is straightforward to check that whenever $\Re\lambda>0$
and $|\Im\lambda|> \tilde a|\lambda|$, then
$\tilde a|\mu|\le |\lambda-\mu|$. Then \qref{Sest} follows,
proving that $B$ is sectorial.

That $(\lambda-B)^{-1}$ is compact for $\lambda\notin\sigma(B)\cup\sigma(B_0)$
follows from the compactness of $(\lambda-B_0)^{-1}$ together
with the identity
\[
(\lambda-B)^{-1} =
(\lambda-B_0)^{-1} + (\lambda-B_0)^{-1} B_1(\lambda-B)^{-1}  .
\]
It follows that the spectrum of $B$ is discrete, consisting only of
isolated eigenvalues of finite multiplicity \cite[III.6.29]{Ka}.

Suppose now that $(\lambda-B)\uu=0$ for some non-zero $\uu\in D(B)$,
so $\lambda\uu=-\PP\lap\uu-\grad\div\uu.$
Then the function $w=\div\uu$ satisfies $\lambda w=-\lap w$ in
$\Omg$, $\nder w=0$ on $\Gam$, i.e., $(\lambda+A)w=0$ (see \qref{domA}).
So if $\lambda\notin \sigma(-A)\subset \R_+$, then $\div\uu=0$, and since
$\ndot\uu=0$ on $\Gam$ we have $\uu=\PP\uu$. Then
\[
\lambda\<\uu,\uu\> =
\<-\PP\lap\uu,\uu\>=
\<-\lap\uu,\PP\uu\>=
\|\grad\uu\|^2 ,
\]
so $\lambda>0$. If $\lambda=0$, then $\div\uu$ is constant, but
$\int_\Omg\div\uu=0$ so $\div\uu=0$ and arguing as above we infer $\uu=0$.
Hence $0$ is not an eigenvalue, and so $0$ is in the resolvent set of $B$.

Lastly, for any $\alpha\ge0$,
given $\uuin\in D(B^\alpha)$, the regularity results for
$e^{-Bt}\uuin$ are standard consequences of the fact that $B^\alpha$
is an isomorphism between its domain and $X$ and commutes with $e^{-Bt}$
\cite[p. 74, Theorem 6.13]{Pa}. For uniqueness, see \cite{B}.
$\square$

\noind {\bf Remark 3.} 
The equation $B\uu=\ff$ has an interesting interpretation in terms
of a stationary Stokes system with prescribed divergence. 
Given any $\ff\in L^2(\Omg,\R^N)$ there is a unique $\uu\in\uspace$
such that $B\uu=\ff$, since $0$ is in the resolvent set of $B$ by the
above theorem.  We can write $\P\ff=\ff+\grad q$ where $q\in
H^1(\Omg)$ with $\int_\Omg q=0$.
Since $B\uu=-\P\lap\uu-\grad\div\uu$, we have $\P\ff=-\P\lap\uu$,
so $\grad(q+\div\uu)=0$. 
Let $\ps$ be the Stokes pressure associated with $\uu$.
Then $(\uu,\ps)$ form a solution to the Stokes system
  \begin{align}
  - \lap \uu + \grad \ps = \ff \quad \text{ in } \Omg,  \label{ste.stokes.1} \\
  - \div \uu = q \quad \text{ in } \Omg, \label{ste.stokes.2} \\
  \uu =0 \quad \text{ on } \Gam. \label{ste.stokes.3}
  \end{align}
As a corollary, we can characterize the domains of positive
integer powers of $B$ by
using the regularity theory for the stationary Stokes
equation (see for example \cite[p. 123, theorem 1.5.3]{Soh} or
\cite[p. 23, proposition 2.2]{Te}).
\begin{corollary} \label{C.Bm}
  Let $\Omg$ be a bounded domain with $C^{2m}$ boundary $\Gam$,
     where $m > 1$ is an integer. Then
  \[ D(B^m) = \{ \uu \; | \; \uu \in H^{2m}(\Omg,\R^N), \; \uu = B \uu =
     \ldots = B^{m-1} \uu =0 \text{ on } \Gam\}.
  \]
\end{corollary}
\noind {\bf Proof}: When $m=1$, the
conclusion is true. Suppose it is true when $m=k-1$. When $m=k$,
take any $\uu \in D(B^k)$. By the definition of $D(B^k)$, 
%which are \[ D(B^k) = \{ \uu \in D(B^{k-1}) \mid B \uu \in
%   D(B^{k-1}) \}, \]
we have $\uu\in D(B^{k-1})$ and $B \uu \in D(B^{k-1})$. 
By assumption, $\ff:=B\uu \in H^{2k-2}$ and 
$ B^{k-1} \uu = B^{k-2}(B\uu) = 0$ on $\Gam$. 
Since $\P$ is bounded on $H^{2k-2}$ \cite[I, Remark 1.6]{Te}
we find that $q\in H^{2k-1}(\Omg)$. 
Now \qref{ste.stokes.1}-\qref{ste.stokes.3} hold,
and we can use the regularity theory of the stationary Stokes
equation cited above to conclude $\uu \in H^{2k}(\Omg,\R^N)$. 
This finishes the proof. $\square$

\noind{\bf Remark 4.} We note that $B$ and $B_0$ have the same domain
and that $D(B_0^{1/2})$ is the closure of $D(B_0)=\uspace$
in norm equivalent to
\[
\|\uu\|^2_{X^{1/2}} = \|B_0^{1/2}\uu\|^2 =
\<-\lap\uu,\uu\>=\|\grad\uu\|^2 ,
\]
the ordinary $H^1$ norm. So $D(B_0^{1/2})=H^1_0(\Omg,\R^N)$.
It is known that if $B$ has {\em bounded imaginary powers} then
for $0<\alpha<1$, $D(B^\alpha)$ can be obtained by interpolation
between $X$ and $D(B)=D(B_0)$ and so $D(B^\alpha)=D(B_0^\alpha)$.
The result that indeed $B+cI$ has bounded imaginary powers for some
$c>0$ apparently
follows from a recent analysis of Abels \cite{Ab} related to the
formulation of Grubb and Solonnikov (although the final result in
\cite{Ab} is stated in terms of the constrained Stokes operator in
divergence-free spaces).

\section{Non-homogeneous side conditions} \label{S.nonhom}
Looking back at the Stokes pressure $\ps$ associated with $\uu$,
one recognizes that the no-slip boundary condition for $\uu$ was
essential for getting the crucial equalities
\qref{psupa}-\qref{stokespair-2} using Lemma~\ref{L.uH2}.
So the important question arises, if general boundary conditions $\uu = \gg$
on $\Gam$ are imposed, do we still have an unconstrained formulation
like \qref{newNSE1}-\qref{newNSE3}? Moreover, what can we say if
the velocity field is not divergence free but is specified as
$\div \uu = h$? Such issues are likely to be relevant in the 
analysis of problems involving complex fluids and low Mach number
flows, for example.

In this section we develop and study an unconstrained formulation
for such non-homogeneous problems.
In this new formulation, $\div \uu -h$ satisfies the
heat equation with no-flux boundary conditions.
The main theorem of this section establishes existence and uniqueness
for strong solutions.

\subsection{An unconstrained formulation}
Consider the Navier-Stokes equations with non-homogeneous boundary
conditions and divergence constraint:
 \begin{align}
  \pa_t \uu + \uu \cdot\! \grad \uu + \nabla p  =  \nu \Delta \uu
  + \ff  & \qquad (t > 0, x \in \Omg ),\label{nh-oldNSE1} \\
  \nabla \cdot \uu  =  h  & \qquad (t \geq 0, x \in \Omg ),
  \label{nh-oldNSE2} \\
  \uu   =  \gg  & \qquad (t \geq 0, x \in \Gam ), \label{nh-oldNSE3} \\
  \uu  = \uuin  & \qquad (t = 0, x \in \Omg ). \label{nh-oldNSE4}
\end{align}
What we have done before can be viewed as replacing the
divergence constraint \qref{nh-oldNSE2} by decomposing the pressure
via the formulae in \qref{EulerP} and \qref{StokesP} in such a way
that the divergence constraint is enforced automatically.
It turns out that in the non-homogeneous case a very similar procedure
works. One can simply use the Helmholtz decomposition to identify
Euler and Stokes pressure terms {\it exactly as before} via the
formulae \qref{EulerP} and \qref{StokesP}, but in addition another
term is needed in the total pressure to deal with the inhomogeneities.
Equation \qref{NSE5} is replaced by
\begin{equation}\label{E.alt2}
 \pa_t \uu + \PP( \udotgrad \uu - \ff - \nu \Delta \uu ) +\grad\pgh
 = \nu \nabla(\div \uu).
\end{equation}
The equation that determines the inhomogeneous pressure $\pgh$ can be
found by dotting with $\grad \phi$ for $\phi\in H^1(\Omg)$,
formally integrating by parts and plugging in the side conditions:
We require
\begin{equation}\label{pgh1}
\<\grad\pgh,\grad \phi\> = -\<\pa_t (\nn \cdot \gg),\phi\>_\Gamma
+\<\pa_t h,\phi\> + \<\nu\grad h,\grad \phi\>
\end{equation}
for all $\phi\in H^1(\Omg)$. With this definition, we see from
\qref{E.alt2} that
% $\div\uu - h$ satisfies a weak form of the heat
% equation with no-flux boundary conditions, namely,
\begin{equation}\label{wk-heat2}
\<\pa_t\uu,\grad \phi\>-\<\pa_t(\nn\cdot \gg),
\phi\>_\Gamma+\<\pa_th,\phi\> =\<\nu\grad(\div\uu-h),\grad \phi\>
\end{equation}
for every $\phi\in H^1(\Omg)$.
This will mean $w:=\div \uu -h $ is a weak solution of
\begin{equation} \label{wk-heat4}
  \pa_t w  = \nu \lap w \; \text{  in } \Omg, \qquad
  \nn \cdot\! \grad w  =  0 \; \text{  on } \Gam,
\end{equation}
with initial condition $w=\div \uuin -h\big|_{t=0}$. So the
divergence constraint will be enforced through exponential
diffusive decay as before (see \qref{nh-diss.iden} below).

The total pressure in \qref{nh-oldNSE1} now has the representation
 \begin{equation} \label{nh-ptotal}
   p = \pe + \nu \ps + \pgh,
 \end{equation}
where the Euler pressure $\pe$ and the Stokes pressure $\ps$ are
determined exactly by \qref{EulerP} and \qref{StokesP} as before, and
$\pgh$ is determined up to a constant by the forcing functions $g$ and
$h$ through the weak-form pressure Poisson equation \qref{pgh1}.
(See Lemma~\ref{L.tuup0} below.) Our unconstrained formulation of
\qref{nh-oldNSE1}-\qref{nh-oldNSE4} then takes the form
 \begin{align}
 \pa_t \uu + \udotgrad \uu +  \nabla \pe +\nu \nabla \ps+ \grad
 \pgh = \nu \Delta \uu + \ff
  &\qquad (t>0,\ x\in\Omg), \label{nh-newNSE1} \\
  \uu = \gg
  &\qquad (t\ge0,\ x\in\Gam), \label{nh-newNSE2}\\
  \uu=\uuin
  &\qquad (t=0,\ x\in\Omg). \label{nh-newNSE3}
 \end{align}

Although the definition of Stokes pressure does not require a no-slip
velocity field, clearly the analysis that we performed in section 2
does rely in crucial ways on no-slip boundary conditions.
So in order to analyze the new unconstrained formulation, we will
decompose the velocity field $\uu$ in two parts.
We introduce a fixed field $\tuu$ in $\Omg\times[0,T]$
that satisfies $\tuu = \gg$ on $\Gam$, and let
\begin{equation} \label{nh-vv}
\vv = \uu - \tuu.
\end{equation}
Then $\vv=0$ on $\Gamma$. With this $\vv$, similar to
\qref{EulerP} and \qref{StokesP} we introduce
\begin{align}
  \grad \qe = (\PP -I) (\vv \cdot\! \grad \vv - \ff), \qquad
  \grad \qs = (I - \PP) \lap \vv - \grad \div \vv. \label{E.qes}
\end{align}
Then we can rewrite \qref{nh-newNSE1} as an equation for $\vv$:
\begin{equation} \label{E.vv2}
   \pa_t \vv + \vv \cdot \! \grad \vv
   + \grad \qe +
    \nu \grad \qs +
   \PP(\tuu \cdot \! \grad \vv + \vv \cdot \! \grad \tuu )
   = \nu \lap \vv + \ff -\tff,
\end{equation}
where
\begin{equation} \label{E.tff}
  \tff :=  \pa_t \tuu
 + \PP(\tuu \cdot\! \grad \tuu -\nu\lap\tuu)
 - \nu\grad\div\tuu+ \grad \pgh.
\end{equation}

\subsection{Existence, uniqueness and dissipation identity}
We will first answer questions concerning the existence and regularity
of $\tuu$ and $\pgh$, then state an existence and uniqueness result
for strong solutions of the unconstrained formulation
\qref{nh-newNSE1}--\qref{nh-newNSE3}.
Let $\Omg$ be a bounded, connected domain in $\R^N$ ($N=2$ or $3$)
with boundary $\Gam$ of class $C^3$.
We assume
%Regarding the forcing functions $\ff$, $\gg$ and $h$ in
%\qref{nh-newNSE1}-\qref{nh-newNSE3} and \qref{pgh1}, we assume
\begin{align}
 \uuin &\in H_{uin} := H^1(\Omg,\R^N) ,  \label{nh-conduin} \\
 \ff &\in H_f := L^2(0,T;L^2(\Omg,\R^N)) , \label{nh-condf} \\
 \gg &\in H_g := H^{3/4}(0,T;L^2(\Gam,\R^N)) \cap L^2(0,T;
H^{3/2}(\Gam,\R^N)) \nonumber \\
& \hspace{1.7cm} \cap \{ \gg \; \big| \; \pa_t (\nn \cdot \gg)
\in L^2(0,T;H^{-1/2}(\Gam)) \}  , \label{nh-condg} \\
 h &\in H_h := L^2(0,T;H^1(\Omg))   \cap H^1(0,T;(H^1)'(\Omg)) .
\label{nh-condh}
 \end{align}
Here $(H^1)'$ is the space dual to $H^1$. We also make the
compatibility assumptions
\begin{align}
&   \gg=\uuin \quad \text{ when $t=0$, $x \in \Gam$},
      \label{nh-compat} \\
&   \<  \pa_t (\nn \cdot \gg) , 1 \>_{\Gam}=\< \pa_t h, 1
\>_{\Omg}.
      \label{nh-conddivu}
\end{align}
We remark that most of the literature on nonhomogeneous Navier-Stokes 
problems \cite{La,Sol,Gr1,GS1,GS2} treats the constrained case with
$h=0$ in $\Omg$ and imposes the condition $\ndot\gg=0$ on $\Gam$.
Amann recently studied very weak solutions without imposing the latter
condition, but only in spaces of very low regularity that
exclude the case considered here \cite{Am2}.

We define
\begin{align}
& V := L^2(0,T;H^2(\Omg,\R^N))\cap H^1(0,T;L^2(\Omg,\R^N)) ,
%\emb C([0,T],H^1(\Omg,\R^N)),
\label{Vdef2} %\\
%& W := L^2(0,T;H^1(\Omg,\R^N))\cap H^1(0,T;(H^1)'(\Omg,\R^N)),
%\emb C([0,T],L^2(\Omg,\R^N))
%\label{Wdef}
\end{align}
and note we have the embeddings 
(\cite[p.~42]{Ta1}, \cite[p.~288]{Ev}, \cite[p.~176]{Te})
\begin{equation}\label{VWemb}
V\emb C([0,T],H^1(\Omg,\R^N)),\quad H_h\emb C([0,T],L^2(\Omg)).
\end{equation}
Notice that we have always used an arrow or tilde to denote a vector. So,
without confusion, we can use $Y(\Omg)$ to denote $Y(\Omg,\R^N)$
or $Y(\Omg)$ as appropriate, and further use $X(Y(\Omg))$ to
denote $X(0,T;Y(\Omg))$.

\begin{lemma} \label{L.tuup0}
Assume \qref{nh-conduin}-\qref{nh-conddivu}. Then, there exists
some $\tuu \in V$ that satisfies
 \begin{equation}
    \tuu(0) = \uuin, \quad  \tuu \big|_{\Gam} = \gg,
       \label{nh-tuu1}
 \end{equation}
and there exists $\pgh \in L^2(H^1(\Omg)/\R)$ satisfying
\qref{pgh1}. Moreover,
 \begin{align}
    \|\tuu\|^2_{V} %_{L^2(H^2(\Omg))}  +  \|\tuu\|^2_{H^1(L^2(\Omg))}
    & \leq C  \big(
    \|\gg\|^2_{H^{3/4}(L^2(\Gam)) \cap H^{3/2}(\Gam)) } +
    % } + \|\gg\|^2_{ L^2( H^{3/2}(\Gam)) } +
    \|\uuin\|^2_{H^1(\Omg)}\big), \label{nh-tuu2} \\
    \|\pgh \|_{L^2(H^1(\Omg)/\R)}  & \leq C
    \( \|\pa_t( \nn \cdot \gg )\|_{L^2(H^{-1/2}(\Gam))}
    + \|h\|_{L^2(H^1)\cap H^1((H^1)')   }\) .
    %_{L^2(H^1)} + \|h\|_{H^1((H^1)')} \).
    \label{nh-prbounds}
 \end{align}
\end{lemma}

\noind {\bf Proof:} (i)
By a trace theorem of Lions and Magenes \cite[vol II, Theorem 2.3]{LM},
the fact $\gg \in H^{3/4}(L^2(\Gam)) \cap L^2(H^{3/2}(\Gam))$
together with \qref{nh-conduin} and
the compatibility condition \qref{nh-compat}
implies the existence of $\tuu \in V$ satisfying \qref{nh-tuu1}.

 (ii) One applies the Lax-Milgram lemma for a.e. $t$ to \qref{pgh1}
 in the space of functions in $H^1(\Omg)$ with zero average.
 We omit the standard details.  $\square$
%Look at \qref{pgh1}. For any $q \in H^1(\Omg),$
% \begin{align}
%  \big|-\< \pa_t (\nn \cdot \gg), q\>_{\Gam} + \< \pa_t h, q\> \big|
%  & \leq \| \pa_t (\nn \cdot \gg)
%    \|_{H^{-1/2}(\Gam)} \|q \|_{H^{1/2}(\Gam)}
%    + \|\pa_t h\|_{(H^1)'(\Omg)}\|q\|_{H^1(\Omg)} \nonumber \\
%  & \leq C \( \| \pa_t (\nn \cdot \gg) \|_{H^{-1/2}(\Gam)}
%    + \|\pa_t h\|_{(H^1)'(\Omg)} \) \| q\|_{H^1(\Omg)} \label{E.vv20}.
% \end{align}
% Then, because of compatibility assumption \qref{nh-conddivu}, we can replace the
% $q$ above by $q-C$ for any
% constant $C$. Since $\Omg$ is connected, we can use the Poincar\'e inequality to replace
% the $\|q\|_{H^1(\Omg)} $ in \qref{E.vv20} by $C \|\grad q\|$.
%
% Going back to \qref{pgh1}, we see that the Lax-Milgram lemma
% implies the existence of
% a $\pgh \in H^1(\Omg)/\R$ for a.e. $t \in [0,T]$.
% Moreover, taking $q = \pgh$, we get
% \begin{equation} \label{nh-prest}
% \|\grad \pgh\| \leq C \( \|\pa_t (\nn \cdot \gg)\|_{H^{-1/2}(\Gam)} + \|\pa_t h\|_{(H^1)'(\Omg)}
% \) + \nu \|\grad h\|
% \end{equation}
%Note that the right hand side is in $L^2([0,T]).$ So, $\pgh \in
%L^2(H^1(\Omg)/\R).$

\begin{theorem}\label{T.nonhomog}
Let $\Omg$ be a bounded, connected domain in $\R^N$ ($N=2$ or $3)$ and
assume \qref{nh-conduin}-\qref{nh-conddivu}. Then
there exists $T^*>0$ so that a unique strong solution of
\qref{nh-newNSE1}-\qref{nh-newNSE3} exists on $[0,T^*]$, with
\begin{align*}
\uu\in&
L^2(0,T^*;H^2(\Omg,\R^N))\cap H^1(0,T^*;L^2(\Omg,\R^N)),\\
p=&\nu\ps+ \pe + \pgh \in L^2(0,T^*;H^1(\Omg)/\R) ,
\end{align*}
where $\pe$ and $\ps$ are defined in \qref{EulerP} and
\qref{StokesP} after introducing the $\tuu$ and $\pgh$ from
Lemma~\ref{L.tuup0}. Moreover, $\uu\in C([0,T^*],H^1(\Omg,\R^N))$
and
$$ \div \uu-h \in L^2(0,T^*;H^1(\Omg)) \cap
H^1(0,T^*;(H^1)'(\Omg))
$$
is a smooth solution of
the heat equation for $t>0$ with no-flux boundary conditions.
The map $t\mapsto\|\div\uu-h\|^2$ is smooth for $t>0$
and we have the dissipation identity
\begin{equation} \label{nh-diss.iden}
  \frac{d}{dt}\frac12\|\div \uu-h \|^2 +  \nu \|\grad (\div \uu -h)\|^2 =0.
\end{equation}
If we further assume
$h \in H_{h.s}:=L^2(0,T;H^2(\Omg))\cap H^1(0,T;L^2(\Omg))$ % V % L^2(0,T^*;H^2(\Omg)) \cap H^1(0,T^*;L^2(\Omg))$
and $\div \uuin \in H^1(\Omg)$,
then
\[
\div \uu \in L^2(0,T^*;H^2(\Omg)) \cap H^1(0,T^*;L^2(\Omg)).
\]
\end{theorem}
\bigskip

\noind {\bf Proof:} First rewrite
\qref{nh-newNSE1} as \qref{E.vv2}. Then we note that there are
only two differences between \qref{E.vv2} and \qref{newNSE1}:

(i) There is an extra forcing term $\tff$ in
\qref{E.vv2}. But by Lemma~\ref{L.tuup0}, all terms in $\tff$
are known to be in $L^2(L^2(\Omg))$ and thus they won't be a problem.

(ii) Equation \qref{E.vv2} has some extra linear terms:
\begin{equation}\label{newterms}
  \PP(\tuu \cdot\! \grad \vv + \vv \cdot\! \grad \tuu ).
\end{equation}
We know $\tuu \in V \emb C([0,T],H^1(\Omg,\R^N))$,
so we can discretize these terms explicitly by setting
$\tuu^n=\tuu(n\Dt)$.
Similar to \qref{E.nlin}, we get
\begin{align}
  \| \PP(\tuu \cdot\! \grad \vv  ) \|^2 \le \eps \| \lap \vv \|^2 +
  \frac{C}{\eps} \|\tuu \|^4_{H^1} \|\grad \vv \|^2. \label{term.1}
\end{align}
We estimate the other term in \qref{newterms} by using
Gagliardo-Nirenberg inequalities \cite[Thm.\ 10.1]{Fr} and the Sobolev
embeddings of $H^1$ into $L^3$ and $L^6$:
\begin{equation}
\|\vv\|_{L^\infty} \leq
\begin{cases}
C \| \lap \vv \|_{L^{3/2}}^{1/2} \|\vv \|_{L^3}^{1/2}
 \leq C \| \lap \vv \|^{1/2} \| \grad \vv \|^{1/2} & (N=2),
 \\[4pt]
C \| \lap \vv \|^{1/2} \|\vv \|_{L^6}^{1/2}
 \leq C \| \lap \vv \|^{1/2} \| \grad \vv \|^{1/2} & (N=3).
\end{cases}
\end{equation}
Then for $N=2$ and $3$ we have
\begin{align}
  \| \PP(\vv \cdot \! \grad \tuu  ) \|^2 \le \|\vv \|^2_{L^\infty}
  \|\grad \tuu\|^2
  \le \eps \| \lap \vv \|^2 +
  \frac{C}{\eps} \| \tuu \|^4_{H^1} \|\grad \vv \|^2. \label{term.2}
\end{align}
With these estimates, the rest of the proof of existence and
uniqueness is essentially the same as that of Theorem~\ref{existence},
and therefore we omit the details.

To prove the regularity of $\div \uu$, we argue in a manner similar
to the proof of Theorem~\ref{existence}. We go from
\qref{nh-newNSE1} to \qref{E.alt2} by using \qref{EulerP} and
\qref{StokesP}. Then using \qref{pgh1} we get \qref{wk-heat2}
for any $\phi\in H^1(\Omg)$.
With $w=\div\uu-h$, taking $\phi\in D(A)$ as in \qref{domA}, we have
\begin{equation}
\<w,\phi\> = \<\nn\cdot\gg,\phi\>_\Gam -\<\uu,\grad\phi\>-\<h,\phi\>,
\end{equation}
therefore $t\mapsto\<w,\phi\>$ is absolutely continuous, and
\qref{wk-heat2} yields
$(d/dt)\<w,\phi\> = \<w,A\phi\>$ for a.e. $t$.
This means $w$ is a weak solution in the sense of Ball \cite{B}, and
the rest of the proof goes as before.

If we further assume
$h \in H_{h.s}$ % V$ %L^2(0,T^*;H^2(\Omg)) \cap H^1(0,T^*;L^2(\Omg))$
and $\div \uuin \in H^1(\Omg)$,
then $w(0)\in H^1(\Omg)$. We claim
\begin{equation}\label{halfdom}
H^1(\Omg)=D((-A)^{1/2}).
\end{equation}
Then semigroup theory yields $w\in C([0,T^*],D((-A)^{1/2}))$,
so since
\begin{equation}
0=\<-\lap w,\pa_t w-\nu\lap w\> = \frac{d}{dt}\frac12 \|\grad w\|^2
+\nu\|\lap w\|^2
\end{equation}
for $t>0$,
we deduce $w\in L^2(0,T^*;H^2(\Omg)) \cap H^1(0,T^*;L^2(\Omg))$,
and $\div\uu$ is in the same space.

To prove \qref{halfdom}, note $X:=D((-A)^{1/2})$ is the closure of
$D(A)$ from \qref{domA} in the norm given by
\[
\|w\|^2_X = \|w\|^2+\|(-A)^{1/2}w\|^2 =
\<(I-\nu\lap)w,w\> = \int_{\Omg}|w|^2+\nu|\grad w|^2.
\]
Clearly $X\subset H^1(\Omg)$. For the other direction, let $w\in
H^1(\Omg)$ be arbitrary. We may suppose $w\in C^\infty(\bar\Omg)$ since
this space is dense in $H^1(\Omg)$. Now we only need to construct a
sequence of $C^2$ functions $w_n\to0$ in $H^1$ norm with
$\nder w_n=\nder w$ on $\Gam$.  This is easily accomplished
using functions of the form $w_n(x)= \xi_n(\dist(x,\Gam))\nder w(x)$,
where $\xi_n(s)=\xi(ns)/n$ with $\xi$ smooth and satisfying
$\xi(0)=0$, $\xi'(0)=1$ and $\xi(s)=0$ for $s>1$.
This proves \qref{halfdom}.

We can prove the uniqueness
by the same method as in Theorem~\ref{existence}.
$\square$

\section{Isomorphism theorems for non-homogeneous Stokes systems}
\label{S.stokes}

Drop the nonlinear term and consider the non-homogeneous Stokes system:
   \begin{align}
     \pa_t \uu + \nabla p  - \nu \Delta \uu = \ff &
\qquad (t>0,\ x\in \Omg), \label{E.stokes1}\\
     \nabla \cdot \uu  = h &
\qquad (t \geq 0,\ x\in \Omg), \label{E.stokes2} \\
     \uu  = \gg &
\qquad (t \geq 0,\ x\in \Gam), \label{E.stokes3}\\
     \uu  = \uuin &
\qquad (t=0,\ x\in \Omg). \label{E.stokes4}
   \end{align}
The unconstrained formulation is
 \begin{align}
  \pa_t \uu +\grad p
%  + (I-\PP)\ff + \nu \nabla \ps+ \grad \pgh
 - \nu \Delta \uu = \ff
   &\qquad (t>0,\ x\in\Omg), \label{nh-newStokes1} \\
  \uu = \gg
  &\qquad (t\ge0,\ x\in\Gam), \label{nh-newStokes2}\\
  \uu=\uuin
  &\qquad (t=0,\ x\in\Omg), \label{nh-newStokes3}
 \end{align}
with
 \begin{equation}\label{nh-pdef}
 \grad p = (I-\PP)\ff + \nu \nabla \ps+ \grad \pgh ,
 \end{equation}
 where $\ps$ and $\pgh$ are defined as before
 via \qref{StokesP} and \qref{pgh1}.

The aim of this section is to obtain an isomorphism between the
space of solutions and the space of data $\{\ff,\gg,h,\uuin \}$,
for this unconstrained formulation and for the original Stokes system.
In examining this question we are motivated by the classic works of Lions and
Magenes~\cite{LM} which provide a satisfactory description of the
correspondence between solutions and data for elliptic boundary value
problems.  In the spirit of these results, a satisfactory theory of a
given system of partial differential equations should describe exactly
how, in the space of all functions involved, the manifold of solutions
can be parametrized.  Yet we are not aware of any such complete
treatment of the non-homogeneous Stokes system. (See further remarks
on this issue below.)

First we consider the mapping from data to solution.
Thanks to the absence of the nonlinear term, we can repeat
much easier what we did in the proof of
Theorems~\ref{existence} and \ref{T.nonhomog} and get the global
existence and uniqueness of a strong solution
of \qref{nh-newStokes1}-\qref{nh-pdef} under the same assumptions
as Theorem~\ref{T.nonhomog}.
The data $\{\ff,\gg,h,\uuin \}$ lie inside the space
\begin{equation}
\Pi_{F}:=
 H_f \times H_g \times H_{h} \times H_{uin}
\end{equation}
from \qref{nh-conduin}--\qref{nh-condh}, and
need to satisfy the compatibility conditions
\qref{nh-compat}--\qref{nh-conddivu}.
Corresponding to such data, we get a unique solution $\uu$
of \qref{nh-newStokes1}-\qref{nh-pdef} in the space
\begin{align}
 H_{u} &:= L^2(0,T;H^2(\Omg,\R^N)) \cap H^1(0,T;L^2(\Omg,\R^N)) \nonumber \\
     & \qquad \cap \{ \uu \mid \pa_t (\nn \cdot \uu)|_\Gam \in
          L^2(0,T;H^{-1/2}(\Gam)) \}. \label{S.Hu}
\end{align}
The total pressure $p$ lies in
\begin{equation}
    H_p :=  L^2(0,T;H^1(\Omg)/\R), \label{S.Hp}
\end{equation}
and the pair $\{\uu,p\}$ satisfies \qref{E.stokes1}, \qref{E.stokes3} and
\qref{E.stokes4}.
As in Theorem~\ref{T.nonhomog}, we can show $w=\div \uu -h$ satisfies a
heat equation with no-flux boundary conditions.
Equation \qref{E.stokes2} says that $w=0$, and this
will hold if and only if $w(0)=0$,
i.e., the following additional compatibility condition holds:
\begin{equation}\label{c.hin}
\div\uuin = h(0).
\end{equation}
For the non-homogeneous Stokes system \qref{E.stokes1}--\qref{E.stokes4},
then, we define the data and solution spaces by
\begin{align}
\Pi_{F.c}&:= \Big\{ \{\ff,\gg,h,\uuin \}\in \Pi_F :
   \text{ \qref{nh-compat}, \qref{nh-conddivu} and \qref{c.hin} hold }
   \Big\},\label{s.fc} \\
\Pi_U &:= H_u\times H_p.\label{s.u}
\end{align}

From what we have said so far, we get a map
$\Pi_{F.c}\to\Pi_U$ by solving the unconstrained system
\qref{nh-newStokes1}--\qref{nh-pdef}.
Due to the absence of nonlinear terms,
the estimates in the proof ensure that this map is bounded.
In the other direction, given $\{\uu,p\}\in \Pi_U$, we simply define
$\{\ff,\gg,h,\uuin \}$ using \qref{E.stokes1}--\qref{E.stokes4}
and check that this lies in $\Pi_{F.c}$.

Note that in Theorem~\ref{T.nonhomog}, one has more regularity on
$\div \uu$ if one assumes more on $\div \uuin$ and $h$.
Correspondingly, like $H_{h.s}$ defined in
Theorem~\ref{T.nonhomog}, we introduce spaces of stronger
regularity by
\begin{align}
   & H_{uin.s}  :=
   H^1(\Omg,\R^N) \cap \{ \uuin \; \big| \; \div \uuin \in H^1(\Omg) \},
   \label{Space.Hs.uin}\\
%   & H_{h.s}  := L^2(0,T;H^2(\Omg)) \cap H^1(0,T;L^2(\Omg)) ,
%            \label{Space.Hsh}\\
& \Pi_{F.s}:= H_f \times H_g \times H_{h.s} \times H_{uin.s}.
\end{align}
The solution $\uu$ then lies in
\begin{align}
& H_{u.s} := L^2(0,T;H^2(\Omg,\R^N)) \cap H^1(0,T;L^2(\Omg,\R^N)) \nonumber \\
     & \hspace{2cm} \cap \{ \uu \; \big| \; \div \uu \in L^2(0,T;H^2(\Omg))
     \cap H^1(0,T;L^2(\Omg)) \}. \label{Space.Hsu}
\end{align}
(Note, if $\uu\in H_{u.s}$ then $\pa_t\uu\in L^2(\Hdiv)$ so
$\nn\cdot\pa_t\uu\in L^2(H^{-1/2}(\Gam))$.)
So as an alternative to the spaces in \qref{s.fc}--\qref{s.u},
we also obtain an isomorphism between the data and solution spaces
with stronger regularity defined by
\begin{align}
\Pi_{F.c.s}&:= \Big\{ \{\ff,\gg,h,\uuin \}\in \Pi_{F.s} :
   \text{ \qref{nh-compat}, \qref{nh-conddivu} and \qref{c.hin} hold }
   \Big\},\label{s.fsc} \\
\Pi_{U.s} &:= H_{u.s}\times H_p.
\label{s.us}
\end{align}

Summarizing, we have proved the following isomorphism theorem for the
non-homogeneous Stokes system \qref{E.stokes1}--\qref{E.stokes4}.
\begin{theorem} \label{T.Stokes}
   Let $\Omg$ be a bounded, connected domain in $\R^N$ with $N$ any positive
   integer $\ge 2$, and let $T>0$. The map
$\{\ff,\gg,h,\uuin \} \mapsto \{\uu,p\}$, given
by solving the unconstrained system \qref{nh-newStokes1}--\qref{nh-pdef},
defines an isomorphism from
$\Pi_{F.c}$ onto $\Pi_U$. The same solution procedure defines
an isomorphism from $\Pi_{F.c.s}$ onto $\Pi_{U.s}$.
\end{theorem}

\medskip \noindent{\bf Remark 5.}
For the standard Stokes system with zero-divergence constraints
$\div \uuin = 0$  and $h=0$,
existence and uniqueness results together with the estimates
\begin{align}
 \sup_{0 \leq t \leq T} & \| \uu (t)\|_{H^1}  +
   \| \uu  \|_{L^2(0,T;H^2)} + \| p \|_{L^2(0,T;H^1/\R)} \nonumber \\
   & \leq C  \( \| \ff \|_{L^2(0,T;L^2)} + \|\uuin\|_{H^1} +
   \|\gg\|_{H^{3/4}(L^2(\Gam))}+ \|\gg\|_{L^2(H^{3/2}(\Gam))}
   \) \label{E.solo}
\end{align}
were obtained in the classic work of Solonnikov
\cite[Theorem 15]{Sol}, where more general $L^p$ estimates were
also proved. (Also see \cite{GS1,GS2}.)
However, instead of the necessary compatibility condition
\begin{equation} \label{E.temp100}
   \int_{\Gam}\ndot \gg=0,
\end{equation}
Solonnikov made the stronger constraining assumption
that both the data $\gg$ and solution $\uu$ have zero normal component
on $\Gam$,
%$$ \nn \cdot \gg =0 \; \text{ on } \Gam,
%$$
and correspondingly his estimates do not contain
a term $\|\pa_t (\nn \cdot
\gg)\|_{L^2(H^{-1/2}(\Gam))}$ on the right hand side of
\qref{E.solo}.
(Note that when $\div \uuin =0$ and $h=0$, we have
$\int_{\Gam}\ndot \gg|_{t=0}=\int_{\Omg} \div \uuin=0$ by
\qref{nh-compat}, whence \qref{E.temp100} is equivalent to \qref{nh-conddivu}.)

\bigskip
\noindent{\bf Remark 6.}
For the unconstrained Stokes system
\qref{nh-newStokes1}-\qref{nh-newStokes3} there is an extra subtlety
in determining an isomorphism from data to solution.
We obtain a unique solution pair $\{\uu,p\}\in\Pi_U$
given any data $\{\ff,\gg,h,\uuin\}\in \Pi_F$ that satisfy only
the compatibility conditions
\qref{nh-compat} and \qref{nh-conddivu} {\em without} \qref{c.hin}.
Consequently the map from data to $\{\uu,p\}$ is not one-to-one.
And, in the other direction, given $\{\uu,p\}$, we can recover
 \begin{equation} \label{E.mapK1}
\ff = \pa_t \uu + \grad p - \nu \Delta \uu,
\qquad \gg = \uu \big|_{\Gam},
\qquad \uuin = \uu|_{t=0}.
 \end{equation}
But how are we to recover $h$? We need to use the fact, that follows
from the definition of $\pgh$ in \qref{pgh1}, that $\div\uu-h$
satisfies a heat equation with no-flux boundary conditions.
In fact, to be able to recover $h$
we need to know one more item, $\hin$, the initial value of
$h$.  We have
\begin{equation} \label{E.mapK2}
  h=\div \uu - w
\end{equation}
where $w$ is the solution of
\begin{equation} \label{E.mapK3}
  \pa_t w  = \nu \lap w \; \text{  in } \Omg, \qquad
  \nn \cdot \grad w  =  0 \; \text{  on } \Gam, \qquad
w(0)=\div \uu|_{t=0} - \hin.
\end{equation}
This procedure
indicates that we should count the triple $\{\uu,p,\hin \}$ as our
solution in order to build an isomorphism with the data. Of course,
the regularity of $\hin$ must match that of $h$,
recalling the embeddings in \qref{VWemb}.

Consequently, we see that solving the unconstrained system
\qref{nh-newStokes1}--\qref{nh-pdef} defines an isomorphism
between the data spaces
\begin{align}
\tilde\Pi_{F.c}&:= \Big\{ \{\ff,\gg,h,\uuin \}\in \Pi_F :
   \text{ \qref{nh-compat} and \qref{nh-conddivu}  hold }
   \Big\},\label{ts.fc} \\
\tilde\Pi_{F.c.s}&:= \Big\{ \{\ff,\gg,h,\uuin \}\in \Pi_{F.s} :
   \text{ \qref{nh-compat} and \qref{nh-conddivu}  hold }
   \Big\},\label{ts.fsc}
\end{align}
and, respectively, the solution spaces for $\{\uu,p,\hin\}$
given by
\begin{align}
 & \Pi_{U.w} = H_{u} \times H_p \times H_{hin} ,
\qquad H_{hin}  =  L^2(\Omg),
     \label{Space.Pi.Uw} \\
 & \Pi_{U.s} = H_{u.s} \times H_p \times H_{hin.s} ,
\qquad H_{hin.s}  =  H^1(\Omg).
     \label{Space.Pi.Us}
\end{align}

\section*{Acknowledgments}
The fact (related to Theorem~\ref{T.N2D})
that $\nder p\in L^2(\Omg_s)$ implies $\grad p\in L^2(\Omg)$
for harmonic $p$ was proved some years ago by Oscar Gonzalez and RLP
(unpublished) through a partitioning and flattening argument.
RLP is grateful for this collaboration.
This material is based upon work supported by the National Science
Foundation under grant no.\ DMS 03-05985 (RLP) and DMS-0107218 (JGL).
JGL and RLP are thankful for the support of the Institute for Mathematical
Sciences at the National University of Singapore.
RLP acknowledges support by the Distinguished Ordway Visitors
Program of the School of Mathematics, and the Institute for
Mathematics and its Applications, at the University of Minnesota.

%%%-------------------------------------

\end{document}